\documentclass[reqno]{amsart}

\usepackage{amssymb,amscd}

        \numberwithin{equation}{section}

\theoremstyle{plain}

\newtheorem{Thm}{Theorem}[section]
\newtheorem{Lem}[Thm]{Lemma}
\newtheorem{Cor}[Thm]{Corollary}
\newtheorem{Prop}[Thm]{Proposition}
\newtheorem{MainThm}{Theorem}
\newtheorem{MainCor}[MainThm]{Corollary}

\theoremstyle{definition}

\newtheorem{Def}[Thm]{Definition}

\theoremstyle{remark}

\newtheorem{Rem}{Remark}
\newtheorem{claim}{Claim}

\newtheorem{ack}{Acknowledgment}

\renewcommand{\sec}[2]{\section{#2}\label{s:#1}}
\newcommand{\ssec}[2]{\subsection{#2}\label{ss:#1}}

\newcommand{\refs}[1]{Section~\ref{s:#1}}
\newcommand{\refss}[1]{Section~\ref{ss:#1}}
\newcommand{\reft}[1]{Theorem~\ref{t:#1}}
\newcommand{\refl}[1]{Lemma~\ref{l:#1}}
\newcommand{\refp}[1]{Proposition~\ref{p:#1}}
\newcommand{\refc}[1]{Corollary~\ref{c:#1}}

\newcommand{\refr}[1]{Remark~\ref{r:#1}}
\newcommand{\refe}[1]{\eqref{e:#1}}

\newcommand{\refmt}[1]{Theorem~\ref{mt:#1}}

\newcommand{\refls}[2]{Lemmas~\ref{l:#1} and~\ref{l:#2}}
\newcommand{\refcs}[2]{Corollaries~\ref{c:#1} and~\ref{c:#2}}
\newcommand{\refmc}[1]{Corollary~\ref{mc:#1}}

\newenvironment{thm}[1]{\begin{Thm} \label{t:#1}}{\end{Thm}}
\renewcommand{\th}[1]{\begin{thm}{#1}}
\renewcommand{\eth}{\end{thm}}

\newenvironment{mainthm}[1]{\begin{MainThm}\label{mt:#1}}{\end{MainThm}}
\newcommand{\mainth}[1]{\begin{mainthm}{#1}}
\newcommand{\emainth}{\end{mainthm}}

\newenvironment{lemma}[1]{\begin{Lem}\label{l:#1}}{\end{Lem}}
\newcommand{\lem}[1]{\begin{lemma}{#1}}
\newcommand{\elem}{\end{lemma}}

\newenvironment{propos}[1]{\begin{Prop}\label{p:#1}}{\end{Prop}}
\newcommand{\prop}[1]{\begin{propos}{#1}}
\newcommand{\eprop}{\end{propos}}

\newenvironment{corol}[1]{\begin{Cor}\label{c:#1}}{\end{Cor}}
\newcommand{\cor}[1]{\begin{corol}{#1}}
\newcommand{\ecor}{\end{corol}}

\newenvironment{maincorol}[1]{\begin{MainCor}\label{mc:#1}}{\end{MainCor}}
\newcommand{\maincor}[1]{\begin{maincorol}{#1}}
\newcommand{\emaincor}{\end{maincorol}}

\newenvironment{defini}[1]{\begin{Def}\label{d:#1}}{\end{Def}}
\newcommand{\defi}[1]{\begin{defini}{#1}}
\newcommand{\edefi}{\end{defini}}

\newenvironment{remark}[1]{\begin{Rem}\label{r:#1}}{\end{Rem}}
\newcommand{\rem}[1]{\begin{remark}{#1}}
\newcommand{\erem}{\end{remark}}

\newcommand{\prf}{ \begin{proof} }
\newcommand{\eprf}{ \end{proof} }

\newcommand{\calA}{\text{$\mathcal{A}$}}
\newcommand{\calB}{\text{$\mathcal{B}$}}

\newcommand{\calF}{\text{$\mathcal{F}$}}

\newcommand{\calH}{\text{$\mathcal{H}$}}

\newcommand{\calZ}{\text{$\mathcal{Z}$}}

\newcommand{\bfb}{{\mathbf b}}		\newcommand{\bfB}{{\mathbf B}}
		
\newcommand{\bfd}{{\mathbf d}}		\newcommand{\bfD}{{\mathbf D}}
\newcommand{\bfe}{{\mathbf e}}		\newcommand{\bfE}{{\mathbf E}}

\newcommand{\bfo}{{\mathbf o}}

\newcommand{\bfz}{{\mathbf z}}		\newcommand{\bfZ}{{\mathbf Z}}

\newcommand{\bfzero}{{\mathbf 0}}

\newcommand{\bfdelta}{\boldsymbol{\delta}}		\newcommand{\bfDelta}{\boldsymbol{\Delta}}

		\newcommand{\bfOmega}{\boldsymbol{\Omega}}

        \newcommand{\field}[1]{\text{$\mathbb{#1}$}}
        
        \newcommand{\Z}{\field{Z}}
        
        \newcommand{\R}{\field{R}}

\newcommand{\im}{\operatorname{im}}
\newcommand{\dom}{\operatorname{dom}}
\newcommand{\cl}{\operatorname{cl}}
\newcommand{\id}{\operatorname{id}}

\newcommand{\cinf}{\text{$C^\infty$}}
\newcommand{\lar}{\text{$\longrightarrow$}}

\newcommand{\norm}[1]{\left\|#1\right\|}

\newcommand{\D}{\text{$\Delta$}}

\begin{document}

\bibliographystyle{plain}

%-------------------------------------------------------------------
%-------------------------------------------------------------------

\title[Adiabatic limits and spectral sequences]{Adiabatic limits and spectral
sequences for Riemannian foliations}
\author[J.A. \'Alvarez L\'opez]{Jes\'us A. \'Alvarez L\'opez}
\address{Departamento de Xeometr\'{\i}a e Topolox\'{\i}a\\
         Facultade de Matem\'aticas\\
         Universidade de Santiago de Compostela\\
         15706 Santiago de Compostela\\ Spain}
\email{alvarez@zmat.usc.es}
\thanks{Partially supported by Xunta de Galicia, grant XUGA20701B97}

\author[Y.A. Kordyukov]{Yuri A. Kordyukov}
\address{Department of Mathematics\\
         Ufa State Aviation Technical University\\
         12~K.~Marx str.\\ 450025 Ufa\\ Russia}
\email{yurikor@math.ugatu.ac.ru}
\thanks{Partially supported by Direcci\'on General de Ense\~nanza Superior e Investigaci\'on Cient\'\i fica
(Spain), sabbatical grant SAB1995-0717}

\begin{abstract}
For general Riemannian foliations, spectral asymptotics of the Laplacian
is studied when the metric on the ambient manifold is blown up in directions normal to the
leaves (adiabatic limit). The number of ``small'' eigenvalues is given in terms of the
differentiable spectral sequence of the foliation. The asymptotics of the corresponding
eigenforms also leads to a Hodge theoretic description of this spectral sequence. This is an
extension of results of Mazzeo-Melrose and R.~Forman.
\end{abstract}

\maketitle

\tableofcontents

%--------------------------------------------------------------
%--------------------------------------------------------------

\sec{intro}{Introduction and main results}

Let $\calF$ be a \cinf~foliation on a closed Riemannian manifold $(M,g)$, and let $T\calF\subset TM$ denote
the subbundle of vectors tangent to the leaves. Then the metric $g$ can be written as an orthogonal sum,
$g=g_\perp\oplus g_F$, with respect to the decomposition $TM=T\calF^\perp\oplus T\calF$; i.e., $g_\perp,g_F$
are the restrictions of $g$ to $T\calF^\perp,T\calF$, respectively. By introducing a parameter $h>0$, we can
define a family of metrics 
\begin{equation}\label{e:gh}
g_h=h^{-2}g_\perp\oplus g_F\;.
\end{equation}
The ``limit'' of the Riemannian manifolds $(M,g_h)$ as $h\downarrow0$ is what
is known as {\em adiabatic limit\/}. Observe that, in a foliation chart, the
plaques get further from each other as $h\downarrow0$. This form of the adiabatic
limit was introduced by E.~Witten in \cite{Witten85} for Riemannian
bundles over the circle. Witten investigated the limit of the eta
invariant of the Dirac operator. This question was also considered in \cite{BismutFreedI},
\cite{BismutFreedII} and
\cite{Cheeger87}, and extended to general Riemannian bundles in
\cite{BismutCheeger} and \cite{Dai91}. 

New properties of adiabatic limits were discovered by Mazzeo and Melrose
for the case of Riemannian bundles, relating them
to the Leray spectral sequence \cite{MazzeoMelrose}. This work was used in
\cite{Dai91}, and further developed by R.~Forman in
\cite{Forman95}, where the very general setting of any pair of complementary
distributions is considered. Nevertheless the most interesting results of
\cite{Forman95} are only proved for foliations satisfying very restrictive
conditions. The ideas from \cite{MazzeoMelrose} and \cite{Forman95} were also applied
to the Rumin's complex by Z.~Ge \cite{Ge94}, \cite{Ge95}. 

For a general \cinf\ foliation \calF\ on $M$, the role of (the differentiable version of) the Leray spectral
sequence is played by the so called {\em differentiable spectral sequence\/} $(E_k,d_k)$, which converges to
the de~Rham cohomology of $M$. The definition of $(E_k,d_k)$ is given by filtering the de~Rham complex
$(\Omega,d)$ of $M$ as in the bundle case: A differential form $\omega$ of degree $r$ is said to be of
filtration $\ge k$ if it vanishes whenever $r-k+1$ of the vectors are tangent to the leaves; that is, roughly
speaking, if $\omega$ is of degree $\ge k$ transversely to the leaves. Moreover the
\cinf~topology of $\Omega$ induces a topological
vector space structure on each term $E_k$ such that $d_k$ is continuous. 
A subtle problem here is that $E_k$ may not be Hausdorff
\cite{Haefliger80}. So it makes sense to consider the subcomplex given by the
closure of the trivial subspace,
$\bar0_k\subset E_k$, as well as the quotient complex ${\widehat E}_k=E_k/\bar0_k$, whose
differential operator will be also denoted by $d_k$.

The differentiable spectral sequence is known to satisfy certain good properties for the so
called {\em Riemannian foliations\/}, which are the foliations with ``rigid
transverse dynamics''; i.e., foliations with isometric holonomy for some
Riemannian metric on smooth transversals. A characteristic property of Riemannian foliations is the existence
of a so called {\em bundle-like metric\/} on the ambient manifold, which means that the foliation is
locally defined by Riemannian submersions \cite{Reinhart59}, \cite{Molino82}, \cite{Molino88}. For such
foliations, each term
$E_k$ is Hausdorff of finite dimension if $k\geq2$, and $H(\bar0_1)=0$ \cite{Masa92},
\cite{AlvKordy1}. So 
$E_k\cong{\widehat E}_k$ for $k\geq2$.
Moreover it was recently proved by X.~Masa and the first author that, for $k\geq2$, the terms $E_k$
are homotopy invariants of Riemannian foliations \cite{AlvMasa1}---this generalizes  previous work
showing the topological invariance of the so called {\em basic cohomology\/}
\cite{KacimiNicolau}---. 

Besides the requirement that \calF\ has to be a Riemannian foliation, the mentioned restrictive
hypothesis of R.~Forman in
\cite{Forman95} is that the positive spectrum of the
``leafwise Laplacian'' on $\Omega$ must be bounded away from
zero\footnote{The leafwise Laplacian is what will be denoted by
$\D_0$ in this paper.}. Both conditions together are so strong that the only
examples we know are Riemannian foliations with compact
leaves; i.e., Seifert bundles.
The purpose of our paper is to generalize Forman's work to
arbitrary Riemannian foliations. To state our first main result, let $\D_{g_h}$
denote the Laplacian defined by
$g_h$ on differential forms, and let 
$$0\leq\lambda_0^r(h)\leq\lambda_1^r(h)\leq\lambda_2^r(h)\leq\cdots$$
denote its spectrum on $\Omega^r$, taking multiplicities into
account. It is well known that the eigenvalues of the Laplacian on differential forms vary
continuously under continuous perturbations of the metric \cite{Cheeger83}, and thus the
``branches'' of eigenvalues $\lambda_i^r(h)$ depend continuously on $h>0$. In this paper, we shall only
consider the ``branches'' $\lambda_i^r(h)$ that are convergent to zero as $h\downarrow0$; roughly speaking,
the ``small'' eigenvalues. The asymptotics as $h\downarrow0$ of these metric invariants is related to the
differential invariant $\widehat E_1^r$ and the homotopy invariants $E_k^r$, $k\geq2$, as follows.

\mainth{small eigenvalues} 
With the above notation, for Riemannian foliations on closed Riemannian manifolds we have
\begin{align}
\dim{\widehat E}_1^r&=
\sharp\,\left\{i\ \left|\ \lambda_i^r(h)\in O\left(h^2\right)
\quad\text{as}\quad h\downarrow 0\right.\right\}\;,\label{e:small eigenvalues 1}\\
\dim E_k^r&=
\sharp\,\left\{i\ \left|\ \lambda_i^r(h)\in O\left(h^{2k}\right)
\quad\text{as}\quad h\downarrow 0\right.\right\}\;,\quad k\geq2\;.
\label{e:small eigenvalues k}
\end{align}
\emainth

As a part of the proof of \refmt{small eigenvalues}, and also because of its own
interest, we shall also study the asymptotics of eigenforms of $\D_{g_h}$ corresponding to ``small''
eigenvalues. This study was begun in \cite{MazzeoMelrose} for the case of Riemannian bundles, and continued
in \cite{Forman95} for general complementary distributions. From both \cite{MazzeoMelrose} and
\cite{Forman95}, certain rescaling $\Theta_h$ of differential forms, depending on $h>0$, is crucial to study
this asymptotics. 

The following well known technicality will be useful to explain $\Theta_h$. The
decomposition $TM=T\calF^\perp\oplus T\calF$ induces a bigrading
\begin{equation}\label{e:bigrading}
\bigwedge TM^\ast
=\bigoplus_{u,v}
\left(\bigwedge^uT\calF^{\perp\ast}\otimes\bigwedge^vT\calF^\ast\right)\;;
\end{equation}
roughly speaking, here $u$ denotes transverse degree and $v$ tangential degree. Then the bigrading of
$\Omega$ is defined by considering \cinf\ sections of \refe{bigrading}; i.e., each $\Omega^{u,v}$ is
the space of \cinf\ sections of $\bigwedge^uT\calF^{\perp\ast}\otimes\bigwedge^vT\calF^\ast$. Then the
de~Rham derivative and coderivative decompose as sum of bihomogeneous components,
\begin{equation}\label{e:d,delta}
d=d_{0,1}+d_{1,0}+d_{2,-1}\;,\quad \delta=\delta_{0,-1}+\delta_{-1,0}+\delta_{-2,1}\;,
\end{equation}
where the double subindex denotes the corresponding bidegree (see e.g. \cite{Alv1}); observe that
$d_{i,j}^\ast=\delta_{-i,-j}$. 

Now define $\Theta_h\omega=h^u\omega$ if $\omega\in\bigwedge TM^\ast$ is of transverse degree $u$. 
As pointed out in \cite{MazzeoMelrose} and \cite{Forman95}, such a $\Theta_h$ is an isometry of Riemannian
vector bundles $\left(\bigwedge TM^\ast,g_h\right)\to\left(\bigwedge TM^\ast,g\right)$, where
$g,g_h$ also denote the metrics induced by $g,g_h$ on $\bigwedge TM^\ast$. So we get an isomorphism, also
denoted by $\Theta_h$, between the corresponding Hilbert spaces of $L^2$~sections because the volume 
elements induced by the metrics $g_h$ are multiples of each other. Thus our setting is moved
via $\Theta_h$ to the fixed Hilbert space of square integrable differential forms on $M$ with the inner
product induced by $g$; this Hilbert space is denoted by $\bfOmega$ in this paper. Concretely, we have the
``rescaled derivative'' $d_h=\Theta_hd\Theta_h^{-1}$, whose $g$-adjoint is the ``rescaled coderivative''
$\delta_h=\Theta_h\delta_{g_h}\Theta_h^{-1}$. It is
easy to verify that
\begin{equation}\label{e:dh}
d_h=d_{0,1}+hd_{1,0}+h^2d_{2,-1}
\end{equation}
directly from \refe{d,delta} and the definition of $\Theta_h$. Thus\footnote{Another
way to check \refe{deltah} is by proving directly that
$$\delta_{g_h}=\delta_{0,-1}+h^2\delta_{-1,0}+h^4\delta_{-2,1}\;.$$}
\begin{equation}\label{e:deltah}
\delta_h=\delta_{0,-1}+h\delta_{-1,0}+h^2\delta_{-2,1}\;.
\end{equation}
The ``rescaled Laplacian'' 
$$
\D_h=\Theta_h\D_{g_h}\Theta_h^{-1}=d_h\delta_h+\delta_hd_h
$$ 
is elliptic and essentially self-adjoint in $\bfOmega$. Moreover $\D_h$ has the same
spectrum as
$\D_{g_h}$, and eigenspaces of $\D_{g_h}$ are transformed into eigenspaces of $\D_h$ by $\Theta_h$. We shall
prove that eigenspaces of $\D_h$ corresponding to ``small'' eigenvalues are convergent as
$h\downarrow0$ when the metric $g$ is bundle-like, and the limit is given by a nested sequence of bigraded
subspaces,
$$
\Omega\supset\calH_1\supset\calH_2\supset\calH_3\supset\cdots\supset\calH_\infty\;.
$$
The definition of $\calH_1,\calH_2$ was already given in \cite{AlvKordy1}
as a Hodge theoretic approach to $(E_1,d_1)$ and $(E_2,d_2)$, which is based on our study of leafwise heat
flow. The other spaces $\calH_k$ are defined in this paper as an extension of this Hodge theoretic approach
to the whole spectral sequence $(E_k,d_k)$ (see Sections~\ref{ss:Hodge} and~\ref{ss:Hodge nested seq} for the
precise definition of $\calH_k$). In particular, 
\begin{equation}\label{e:Hodge}
\calH_1\cong\widehat E_1\;,\quad\calH_k\cong E_k\;,\quad k=2,3,\ldots,\infty\;,
\end{equation}
as bigraded topological vector spaces. Thus this sequence stabilizes\footnote{We mean
$\calH_k=\calH_\infty$ for $k$ large enough.} because the differentiable spectral sequence is convergent in a
finite number of steps. The convergence of eigenforms corresponding to ``small'' eigenvalues is precisely
stated in the following result, where $L^2\calH_1$ denotes the closure of $\calH_1$ in $\bfOmega$.

\mainth{asymptotics of eigenforms}
For any Riemannian foliation on a closed manifold with a bundle-like metric, let $\omega_i$ be a sequence in
$\Omega^r$ such that $\|\omega_i\|=1$ and 
\begin{equation}\label{e:asymptotics of eigenforms}
\left\langle\D_{h_i}\omega_i,\omega_i\right\rangle\in o\left(h_i^{2(k-1)}\right)
\end{equation}
for some $k=1,2,3,\dots$ and some sequence $h_i\downarrow0$. Then some subsequence
of the $\omega_i$ is strongly convergent, and its limit is in 
$L^2\calH_1^r$ for $k=1$, and in
$\calH_k^r$ for $k\geq2$.
\emainth

To simplify notation let $m_1^r=\dim\widehat E_1^r$, and let $m_k^r=\dim E_k^r$
for each $k=2,3,\ldots,\infty$. Thus \refmt{small eigenvalues} establishes
$\lambda_i^r(h)\in O\left(h^{2k}\right)$ for $i\leq m_k^r$, yielding $\lambda_i^r(h)\equiv0$ for $i$ 
large enough. For every $h>0$, consider the nested sequence of graded subspaces
$$\Omega\supset\calH_1(h)\supset\calH_2(h)\supset\calH_3(h)\supset\cdots\supset\calH_\infty(h)\;,$$
where $\calH_k^r(h)$ is the space generated by the eigenforms of $\D_h$ corresponding to
eigenvalues $\lambda_i^r(h)$ with $i\leq m_k^r$; in particular, we have
$\calH_k(h)=\calH_\infty(h)=\ker\D_h$ for $k$ large enough. Set also
$\calH_k(0)=\calH_k$. We have $\dim\calH_k^r(h)=m_k^r$ for all $h>0$, so the following result is a
sharpening of \refmt{small eigenvalues}.

\maincor{asymptotics of eigenspaces} 
For any Riemannian foliation on a closed manifold with a bundle-like metric and $k=2,3,\ldots,\infty$, the
assignment $h\mapsto\calH_k^r(h)$ defines a continuous map from $[0,\infty)$ to the space of finite
dimensional linear subspaces of $\bfOmega^r$ for all $r\geq0$. If $\dim\widehat E_1^r<\infty$, then this also
holds for $k=1$.
\emaincor

In \refmc{asymptotics of eigenspaces}, the continuity of $h\mapsto\calH_k^r(h)$ for $h>0$ is a
particular case of the general property that eigenspaces of the Laplacian on closed Riemannian
manifolds vary continuously as subspaces of
$\bfOmega$ when the metric is perturbed $C^0$-continuously \cite{Cheeger83},
\cite{BakerDodziuk97}. On the other hand, the continuity of $h\mapsto\calH_k^r(h)$ at $h=0$ is a
direct consequence of \refmt{asymptotics of eigenforms}.

With an analogous aim, other nested sequences of bigraded subspaces were introduced by Mazzeo-Melrose in
\cite{MazzeoMelrose} and by R.~Forman in \cite{Forman95}, which are respectively denoted by
$$
\Omega\supset{\mathfrak h}_1\supset{\mathfrak h}_2\supset{\mathfrak h}_3\supset\cdots 
\supset{\mathfrak h}_\infty\;,\quad
\Omega\supset{\mathfrak H}_1\supset{\mathfrak H}_2\supset{\mathfrak H}_3\supset\cdots
\supset{\mathfrak H}_\infty
$$
in this paper. These sequences are defined in the following way. According to the expressions \refe{dh} and
\refe{deltah}, we can consider $d_h$ and $\delta_h$ as polynomials on the variable $h$ whose coefficients
are the differential operators $d_{i,j}$ and $\delta_{i,j}$. Thus $d_h$ and $\delta_h$ canonically become
operators on the polynomial algebra $\Omega[h]$, and $\D_h$ as well. Then each ${\mathfrak h}_k$
is the space of differential forms $\omega\in\Omega$ with some extension $\tilde\omega(h)\in\Omega[h]$
satisfying
\begin{equation}\label{e:overlinefrakHk}
\D_h\tilde\omega(h)\in h^k\Omega[h]\;,
\end{equation}
where {\em extension\/} means $\tilde\omega(0)=\omega$.
And each ${\mathfrak H}_k$ is the space of differential forms $\omega\in\Omega$ with some extension
$\tilde\omega(h)\in\Omega[h]$ satisfying
\begin{equation}\label{e:frakHk}
d_h\tilde\omega(h)\in h^k\Omega[h]\;,\quad\delta_h\tilde\omega(h)\in h^k\Omega[h]\;.
\end{equation}
The sequence $\calH_k$ also fits in this kind of description as follows (this is a direct consequence of
\reft{calHk}): Each $\calH_k$ is the space of differential forms $\omega\in\Omega$ having
sequences of extensions $\tilde\omega^1_i(h),\tilde\omega^2_i(h)\in\Omega[h]$ satisfying
\begin{equation}\label{e:calHk}
d_h\tilde\omega^1_i(h)+h^k\Omega[h]\lar0\;,\quad
\delta_h\tilde\omega^2_i(h)+h^k\Omega[h]\lar0
\end{equation}
in $\Omega[h]/h^k\Omega[h]$ as $i\to\infty$. From~\refe{dh},~\refe{deltah},~\refe{frakHk}
and~\refe{calHk} it easily follows that
\begin{gather}
{\mathfrak H}_k\subset{\mathfrak h}_k\subset{\mathfrak H}_{[k/2]}\;.
\label{e:frakHk subset ...}\\
{\mathfrak H}_1=\calH_1\;,\quad{\mathfrak H}_k\subset\calH_k\;,\quad k\geq2\;.
\label{e:frakH1=calH1, ...}
\end{gather}

For the case of Riemannian bundles, Mazzeo and Melrose prove in \cite{MazzeoMelrose} that the sequence
${\mathfrak h}_k$ stabilizes, and ${\mathfrak h}_\infty$ is
the limit of the spaces $\ker\D_h$ as $h\downarrow0$. And for foliations under the restrictive hypothesis
of \cite{Forman95}, R.~Forman proves that the sequence ${\mathfrak H}_k$ is a Hodge
theoretic version of the spectral sequence $(E_k,d_k)$, and describes the limit of the eigenspaces of
$\D_h$ corresponding to ``small'' eigenvalues. This improves the results of Mazzeo-Melrose by \refe{frakHk
subset ...}. But Forman's sequence ${\mathfrak H}_k$ does not have the same important properties for general
Riemannian foliations and bundle-like metrics, as follows from the following result, where the notation
$\calH_k(g)$ and ${\mathfrak H}_k(g)$ is used to emphasize the dependence of
$\calH_k$ and ${\mathfrak H}_k$ on the metric $g$---of course, each $\calH_k(g)$ is independent of $g$ up
to isomorphism by \refe{Hodge}---. 

\mainth{frakHk0,p} 
Let \calF\ be a Riemannian foliation of dimension $p$ on a closed manifold $M$. We have:
\begin{itemize}

\item[$($i$)$] There is a bundle-like metric $g$ on $M$ such that
${\mathfrak H}_2^{0,p}(g)=\calH_2^{0,p}(g)$. 

\item[$($ii$)$] If $\bar0_1^{0,p}\neq0$, then
there is a bundle-like metric $g'$ on $M$ such that ${\mathfrak H}_2^{0,p}(g')=0$. 

\end{itemize}
\emainth

The condition $\bar0_1^{0,p}\neq0$ holds for Kronecker's flows on $T^2$ whose slope is a Liouville's
number \cite{Heitsch75}, \cite{Roger}. This was generalized to linear foliations on tori of arbitrary
dimension in \cite{ArrautdosSantos}. Moreover
$E_2^{0,p}\cong\R$ in these examples \cite{Masa92}, \cite{Alv7}. Therefore \refmt{frakHk0,p} implies that, in
these examples, the dimension of ${\mathfrak H}_2^{0,p}(g)$ changes when appropriately varying the metric
$g$. Thus ${\mathfrak H}_2^{0,p}(g)\not\cong E_2^{0,p}$ for appropriate choices of $g$; that is,
\cite[Corollary~4.4]{Forman95} is not completely right with that generality---the possibility that
$E_1$ may not be Hausdorff is not considered in that paper---.  So far it is rather unknown which
topological or geometric conditions imply $\bar0_1\neq0$ for general Riemannian foliations, but the above
examples suggest that this may happen ``generically''.

A simple argument shows that ${\mathfrak H}_k^r=\calH_k^r$ if $h\mapsto\calH_k^r(h)$ is a
\cinf~map: In this case, any $\omega\in\calH_k^r$ has an extension depending smoothly on $h\geq0$, whose
Taylor polynomial of degree $k$ at zero is easily seen to satisfy \refe{frakHk}, yielding
$\omega\in{\mathfrak H}_k^r$. Therefore, since both $\calH_k$ and ${\mathfrak H}_k$ obviously stabilize at
$k=2$ for flows on surfaces, \refmt{frakHk0,p} shows that the map
$h\mapsto\calH_\infty^1(h)$ is not \cinf\ at $h=0$ for Kronecker's flows on $T^2$ whose slope is a
Liouville's number and appropriate bundle-like metrics. So \cite[Corollary~18]{MazzeoMelrose} and
\cite[Corollary~5.22]{Forman95} have no direct generalizations to arbitrary Riemannian foliations
and bundle-like metrics.

Nevertheless, the arguments of R.~Forman in \cite{Forman95} are
right when $\bar0_1=0$. In particular, Sections~2---4 in \cite{Forman95} show that, in this case, ${\mathfrak
H}_k\cong E_k$ as bigraded vector spaces\footnote{Indeed \cite[Lemma~2.7]{Forman95} is a version of this
isomorphism---it must be pointed out that the notation used in \cite{Forman95} is very different from
ours---.}. Therefore, by~\refe{Hodge} and~\refe{frakH1=calH1, ...}, Forman's arguments prove the following.

\mainth{frakHk} 
Let \calF\ be a Riemannian foliation on a closed manifold $M$. If $\bar0_1=0$, then 
${\mathfrak H}_k(g)=\calH_k(g)$ for every $k\ge1$ and any bundle-like metric $g$ on $M$.
\emainth

\refmt{frakHk0,p}-(ii) is a partial reciprocal of \refmt{frakHk}, and we could conjecture that its
statement holds for any bidegree, but we do not pursue such a result in this paper. A similar question can
be raised about \refmt{frakHk0,p}-(i).

The following are the main ideas of the proofs in this paper. 
The proof of ``$\leq$'' in \refe{small eigenvalues k} (\refmt{small eigenvalues}) has three
main ingredients. The first one is a variational formula for the spectral distribution function of
the Laplacian, which is a consequence of the Hodge decomposition, and was used by Gromov and Shubin
in another setting \cite{GromovShubin}. The second ingredient is a direct sum decomposition that holds for
general spectral sequences---it is kind of an (only linear) Hodge decomposition---. The relation between this
decomposition and the formula of Gromov-Shubin can be easily seen, and leads to the proof. But this can not
be directly applied to the differentiable spectral sequence $(E_k,d_k)$ because of
some technical difficulty (\refr{technical difficulty}). For this reason,
we introduce the third ingredient: The {\em $L^2$~spectral sequence\/}
$(\bfE_k,\bfd_k)$, which is another spectral sequence defined in the very same way as $(E_k,d_k)$ but using
square integrable differential forms. This change of spectral sequence can be made because we show that
$\bfE_k\cong E_k$ for Riemannian foliations and $k\ge2$. The proof of this isomorphism heavily depends on
the Hodge theoretic approach of the terms $E_1$ and $E_2$ that follows from our work \cite{AlvKordy1} on
leafwise heat flow. 

The rest of \refmt{small eigenvalues} is an easy consequence of \refmt{asymptotics of eigenforms}, which in
turn is proved by characterizing the terms $\calH_k$ in the appropriate way to apply certain estimation of
$\D_h$---this estimation is similar to what was done by R.~Forman in \cite{Forman95}---.

\refmt{frakHk0,p} is an easy consequence of the above theorems and other known results about
Riemannian foliations.

Finally, let us mention that a very related study is done in \cite{Kordyukov95b}, where the second author
proves an asymptotical formula for the eigenvalue distribution function of $\D_{g_h}$ in
adiabatic limits for Riemannian foliations. That work establishes relationships with the
spectral theory of leafwise Laplacian and with the noncommutative spectral
geometry of foliations.

\sec{diff spectral seq}{Differentiable spectral sequence}

\ssec{general diff}{General properties}

Let $(\calA,d)$ be a complex with a finite decreasing filtration
$$\calA=\calA_0\supset\calA_1\supset\cdots
\supset\calA_q\supset\calA_{q+1}=0$$
by differential subspaces; i.e. $d(\calA_k)\subset\calA_k$ for all $k$. Recall that the
induced spectral sequence $(E_k,d_k)$ is defined in the
following standard way \cite{McClearly}:
\begin{align*}
Z_k^{u,v}&=\calA_u^{u+v}\cap d^{-1}\left(\calA_{u+k}^{u+v+1}\right)\;,\\
B_k^{u,v}&=\calA_u^{u+v}\cap d\left(\calA_{u-k}^{u+v-1}\right)\;,\\
E_k^{u,v}&=\frac{Z_k^{u,v}}{Z_{k-1}^{u+1,v-1}+B_{k-1}^{u,v}}\;,\\
Z_\infty^{u,v}&=\calA_u^{u+v}\cap \ker d\;,\\
B_\infty^{u,v}&=\calA_u^{u+v}\cap \im d\;,\\
E_\infty^{u,v}&=\frac{Z_\infty^{u,v}}{Z_\infty^{u+1,v-1}+B_\infty^{u,v}}\;.
\end{align*}
In particular
$Z_0^{u,v}=Z_{-1}^{u,v}=\calA_u^{u+v}$. 
We assume $B_{-1}^{u,v}=0$, 
so $E_0^{u,v}=\calA_u^{u+v}/\calA_{u+1}^{u+v}$.
Also, we have $B_u^{u,v}=B_\infty^{u,v}$ and
$Z_{q-u+1}^{u,v}=Z_\infty^{u,v}$ since the filtration of $\calA$ is of length
$q+1$. Each homomorphism $d_k:E_k^{u,v}\to E_k^{u+k,v-k+1}$ is canonically induced
by $d$. 

Now let \calF\ be a \cinf\ foliation
of codimension $q$ on a closed manifold $M$, and $(\Omega,d)$ the de~Rham complex
of $M$. The {\em differentiable spectral sequence\/} $(E_k,d_k)$ of \calF\ is defined by the
decreasing filtration by differential subspaces
$$\Omega=\Omega_0\supset\Omega_1\supset\cdots
\supset\Omega_q\supset\Omega_{q+1}=0\;,$$
where the space of
$r$-forms of filtration degree $\geq k$ is given by
$$\Omega_k^r=\left\{\omega\in\Omega^r\ \left|\ \begin{array}{l}
i_X\omega=0\quad\text{for all}\quad X=X_1\wedge\cdots\wedge X_{r-k+1}\;,\\
\text{where the}\quad X_i\quad\text{are vector
fields tangent}\\ 
\text{to the leaves}\end{array}\right.\right\}\;.$$
Moreover, the \cinf~topology of $\Omega$ canonically induces a topology on each $E_k^{u,v}$,
which becomes a topological vector space. Then each $d_k$ is continuous on
$E_k=\bigoplus_{u,v}E_k^{u,v}$ with the product topology. Thus, for each $k$, we have two new
bigraded complexes: the closure of the trivial subspace $\bar0_k\subset E_k$ and the quotient
$\widehat E_k=E_k/\bar0_k$.

Assume $M$ is endowed with a Riemannian metric, and let $\pi_{u,v}:\Omega\to\Omega^{u,v}$ denote the
induced projection defined by the bigrading of $\Omega$. Define the topological vector spaces
$$z_k^{u,v}=\pi_{u,v}\left(Z_k^{u,v}\right)\;,\quad b_k^{u,v}=\pi_{u,v}\left(B_k^{u,v}\right)\;,\quad
e_k^{u,v}=z_k^{u,v}/b_k^{u,v}\;,\quad e_k=\bigoplus_{u,v}e_k^{u,v}\;.$$
Observe that 
\begin{equation}\label{e:Omega k}
\Omega_k=\bigoplus_{u\geq k}\Omega^{u,\cdot}\;,
\end{equation}
yielding
$$Z_k^{u,v}\cap\ker\pi_{u,v}=Z_{k-1}^{u+1,v-1}\;.$$
Thus the projection
$\pi_{u,v}$ induces a continuous linear isomorphism
$E_k^{u,v}\stackrel{\cong}{\to}e_k^{u,v}$. The operator on $e_k$ that
corresponds to $d_k$ on $E_k$ by the above linear isomorphisms will be denoted by $d_k$ as well. 
We also consider the closure of the trivial subspace, $\bar o_k\subset e_k$,
and the quotient $\hat e_k=e_k/\bar o_k$. We are going to show that $d_k$ is continuous
on $e_k$ for $k=0,1$, and thus $\bar o_k$ and $\hat e_k$ become bigraded complexes in a canonical way.
But, for $k\ge2$, we do not know whether $d_k$ is continuous
on $e_k$, and whether $d_k$ induces differentials on $\bar o_k$ and $\hat e_k$.
This holds at least for Riemannian foliations as easily follows from \reft{H1}-(vii) in \refss{Hodge}.

By comparing bihomogeneous components in the equality $d^2=0$ we get (see e.g. \cite{Alv1}):
\begin{equation}\label{e:d0,12=...}
\left.\begin{array}{c}
d_{0,1}^2=d_{2,-1}^2=d_{0,1}d_{1,0}+d_{1,0}d_{0,1}=0\;,\\[6pt]
d_{1,0}d_{2,-1}+d_{2,-1}d_{1,0}=d_{1,0}^2+d_{0,1}d_{2,-1}+d_{2,-1}d_{0,1}=0\;.
\end{array}\right\}
\end{equation}
The term $d_{2,-1}$ is of order zero, and vanishes if and only if $T\calF^\perp$ is
completely integrable. Moreover from \refe{Omega k} we get
\begin{align}
Z_0^{u,v}&=\Omega_u^{u+v}\;,\label{e:Z0}\\
B_0^{u,v}&=d_{0,1}\left(\Omega^{u,v-1}\right)\oplus
\Omega_{u+1}^{u+v}\;,\label{e:B0}\\ 
Z_1^{u,v}&=\left(\Omega^{u,v}\cap\ker d_{0,1}\right)
\oplus \Omega_{u+1}^{u+v}\label{e:Z1}
\end{align}
as topological vector spaces. So 
\begin{equation}\label{e:z0, b0, z1}
z_0^{u,v}=\Omega^{u,v}\;,\quad
b_0^{u,v}=d_{0,1}\left(\Omega^{u,v-1}\right)\;,\quad
z_1^{u,v}=\Omega^{u,v}\cap\ker d_{0,1}\;,
\end{equation}
and the continuous linear isomorphisms
$E_k^{u,v}\stackrel{\cong}{\to}e_k^{u,v}$, induced by $\pi_{u,v}$, are homeomorphisms too for
$k=0,1$. Thus $\bar0_1\cong\bar o_1$ and $\widehat E_1\cong\hat e_1$ as topological vector spaces,
and $\bar o_1$ and $\hat e_1$ become bigraded complexes with the differential induced by $d_1$.
For this reason, using the spaces $e_1,\bar o_1,\hat e_1$ is rather redundant; we have introduced
these spaces to be compared with the corresponding ones for the
$L^2$~spectral sequence (\refs{L2 spectral seq}), where this does not
obviously hold. Furthermore~\refe{Z0}--\refe{Z1} yield
\begin{equation}\label{e:e0}
(e_0,d_0)=(\Omega,d_{0,1})\;,
\end{equation}
and a canonical isomorphism
\begin{equation}\label{e:e1}
(e_1,d_1)\cong\left(H(\Omega,d_{0,1}),d_{1,0\ast}\right)
\end{equation} 
of topological complexes. Nevertheless we can not go further keeping full control
of the topology. In fact, with this generality, we do not know whether the continuous linear 
isomorphism $E_2^{u,v}\stackrel{\cong}{\to}e_2^{u,v}$, induced by $\pi_{u,v}$, is a
homeomorphism, neither the canonical continuous linear isomorphisms
$E_2\stackrel{\cong}{\to}H(E_1,d_1)$ and $e_2\stackrel{\cong}{\to}H(e_1,d_1)$.

\ssec{Hodge}{Hodge theory of the terms $E_1$ and $E_2$ for Riemannian foliations}

Here, \calF\ is assumed to be a Riemannian foliation and the metric
bundle-like.

The de~Rham coderivative $\delta$ decomposes as sum of
bihomogeneous components $\delta_{i,j}=d_{-i,-j}^\ast$, and the operators
$$
D_0=d_{0,1}+\delta_{0,-1}\;,\quad
\D_0=D_0^2=d_{0,1}\delta_{0,-1}+\delta_{0,-1}d_{0,1}
$$
are essentially self-adjoint in $\bfOmega$ \cite{Chernoff}. But $D_0$ and $\D_0$ are not
elliptic on $M$---avoiding the trivial case where $q=0$---. The closures of $d$, $\delta$, $d_{0,1}$,
$\delta_{0,-1}$, $D_0$ and $\D_0$ in
$\bfOmega$ will be denoted by $\bfd$, $\bfdelta$, $\bfd_{0,1}$, $\bfdelta_{0,-1}$, $\bfD_0$ and
$\bfDelta_0$, respectively.
Then we have the orthogonal decomposition
\begin{equation}\label{e:leafwise Hodge}
\bfOmega=\ker\bfDelta_0\oplus\cl_0\left(\im\bfd_{0,1}\right)
\oplus\cl_0\left(\im\bfdelta_{0,-1}\right)\;,
\end{equation}
where $\cl_0$ denotes closure in $\bfOmega$. Moreover
\begin{gather}
\ker\bfDelta_0=\ker\bfD_0=\ker\bfd_{0,1}\cap\ker\bfdelta_{0,-1}\;,\label{e:ker bfDelta0=ker bfD0}\\
\cl_0\left(\im\bfDelta_0\right)=\cl_0\left(\im\bfD_0\right)
=\cl_0\left(\im\bfd_{0,1}\right)\oplus\cl_0\left(\im\bfdelta_{0,-1}\right)\;\notag
\end{gather}
Thus let $\Pi$, $P$ and $Q$ denote the
orthogonal projections of $\bfOmega$ onto $\ker\bfDelta_0$, $\cl_0\left(\im
\bfd_{0,1}\right)$ and
$\cl_0\left(\im\bfdelta_{0,-1}\right)$, respectively, and set
$\widetilde\Pi=\id-\Pi$, $\widetilde P=\id-P$ and $\widetilde Q=\id-Q$. We
shall also use the notation $W^k\Omega$ for the $k$th Sobolev space completion
of $\Omega$, and let $\cl_k$ denote closure in $W^k\Omega$. Thus $\bfOmega=W^0\Omega$.

\begin{Thm}[{\'Alvarez-Kordyukov \cite{AlvKordy1}}]\label{t:leafwise Hodge}
For each $k\in\Z$, decomposition \refe{leafwise Hodge} restricts to $W^k\Omega$; i.e.,
$$
W^k\Omega=\ker(\D_0\text{ in }W^k\Omega)\oplus\cl_k\left(\im d_{0,1}\right)
\oplus\cl_k\left(\im\delta_{0,-1}\right)
$$
as topological vector space. Thus \refe{leafwise Hodge} also restricts to \cinf\ differential forms; i.e.,
$$
\Omega=\ker\D_0\oplus\overline{\im d_{0,1}}
\oplus\overline{\im\delta_{0,-1}}
$$
with respect to the \cinf~topology, where the bar denotes
\cinf~closure in $\Omega$. In particular $\Pi$, $P$ and $Q$ preserve $\Omega$.
\end{Thm}

From~\refe{e0},~\refe{e1} and \reft{leafwise Hodge}, we get a canonical
isomorphism $\ker\D_0\cong\hat e_1$ of topological vector spaces, induced by the inclusion 
$$
\Omega^{u,v}\cap\ker\D_0\hookrightarrow\Omega^{u,v}\cap\ker d_{0,1}=z_1^{u,v}\;.
$$
So $\ker\D_0\cong\widehat E_1$ as topological vector spaces. As in \cite{AlvKordy1}, let 
\begin{gather*}
{\calH}_1=\ker\D_0=\ker D_0=\ker d_{0,1}\cap\ker\delta_{0,-1}\;,\\
\widetilde{\calH}_1=\overline{\im\D_0}=\overline{\im D_0}
=\overline{\im d_{0,1}}\oplus\overline{\im\delta_{0,-1}}\;,
\end{gather*}
and let $L^2\calH_1=\cl_0\left(\calH_1\right)$ and 
$L^2\widetilde\calH_1=\cl_0\left(\widetilde\calH_1\right)$. From~\refe{ker bfDelta0=ker bfD0} and
\reft{leafwise Hodge} we get
\begin{equation}\label{e:L2calH1}
\ker\bfDelta_0=\ker\bfD_0=L^2\calH_1\;.
\end{equation}
Since
$\D_0$ is bihomogeneous of bidegree $(0,0)$, the bigrading of $\Omega$ restricts to a bigrading of 
${\calH}_1$. Moreover, by \refe{e0}, \refe{e1} and \reft{leafwise Hodge}, the operator $d_1$ on
$\hat e_1$ corresponds to the map $\Pi d_{1,0}$ on
${\calH}_1$, which will be also denoted by $d_1$. Hence 
$H^u({\calH}_1^{\cdot,v},d_1)\cong H^u(\hat e_1^{\cdot,v})\cong H^u({\widehat E}_1^{\cdot,v})$.
Since $\delta_1=\Pi\delta_{-1,0}$ is adjoint of $d_1$ in ${\calH}_1$, the operators $D_1=d_1+\delta_1$ and
$\D_1=D_1^2=d_1\delta_1+\delta_1d_1$ on ${\calH}_1$ are symmetric. Now, let
$\calH_2=\ker\D_1$, which inherits the bigrading from $\Omega$
because $\D_1$ is bihomogeneous of bidegree $(0,0)$.

We also define maps $\tilde d_1$ and $\tilde\delta_1$ on  $\widetilde{\calH}_1$ as follows. First
we define the following bigrading on $\widetilde{\calH}_1$:
$$\widetilde{\calH}_1^{u,v}=\overline{d_{0,1}(\Omega^{u,v-1})}\oplus
\overline{\delta_{0,-1}(\Omega^{u+1,v})}\;.$$ 
Let $\widetilde\Pi_{\cdot,v}$ be the
projection of $\Omega$ onto $\widetilde{\calH}_1^{\cdot,v}$, and set $\tilde
d_1= \widetilde\Pi_{\cdot,v}d$ and $\tilde\delta_1=
\widetilde\Pi_{\cdot,v}\delta$ on
$\widetilde{\calH}_1^{\cdot,v}$, which are adjoint of each other. Consider also 
the symmetric operators $\widetilde{D}_1=\tilde d_1+\tilde\delta_1$ and
$\widetilde{\D}_1=\widetilde{D}_1^2$ on $\widetilde{\calH}_1$. 

The closures of $d_1$, $\delta_1$, $D_1$ and $\D_1$ in $L^2\calH_1$, and of
$\tilde d_1$, $\tilde\delta_1$, $\widetilde{D}_1$ and $\widetilde{\D}_1$ in $L^2\widetilde\calH_1$, will be
respectively denoted by $\bfd_1$, $\bfdelta_1$, $\bfD_1$, $\bfDelta_1$, $\tilde{\bfd}_1$,
$\tilde{\bfdelta}_1$, $\widetilde{\bfD}_1$ and $\widetilde{\bfDelta}_1$.

The following theorem collects the main results of~\cite[Section~7]{AlvKordy1}.

\begin{Thm}[{\'Alvarez-Kordyukov \cite{AlvKordy1}}]\label{t:H1}
We have:
\begin{itemize}

\item[$($i$)$] The operators $D_1$ and $\D_1$ are essentially
self-adjoint in $ L^2\calH_1$, and the operators $\widetilde{D}_1$ and $\widetilde{\D}_1$ are essentially
self-adjoint in $L^2\widetilde\calH_1$. 

\item[$($ii$)$] The spectrums of $\bfD_1$, $\bfDelta_1$, $\widetilde{\bfD}_1$ and
$\widetilde{\bfDelta}_1$ are discrete subsets of \R\ given by eigenvalues of finite multiplicity. 

\item[$($iii$)$] We have the Hodge type decompositions
\begin{align*}
 L^2\calH_1&=\ker\bfDelta_1\oplus\im\bfd_1\oplus\im\bfdelta_1\;,\\
 L^2\widetilde\calH_1&=\im\tilde\bfd_1\oplus\im\tilde{\bfdelta}_1\;,
\end{align*}
as Hilbert spaces with the $L^2$~norm, and moreover
\begin{gather*}
\ker\bfDelta_1=\ker\bfD_1=\ker\bfd_1\cap\ker\bfdelta_1\;,\\
\im\bfDelta_1=\im\bfD_1=\im\bfd_1\oplus\im\bfdelta_1\;,\\
\ker\widetilde{\bfDelta}_1=\ker\widetilde{\bfD}_1=0\;,\quad
\im\bfDelta_1=\im\bfD_1=L^2\widetilde\calH_1\;.
\end{gather*}
Furthermore the operators $\bfDelta_1$ and
$\widetilde{\bfDelta}_1$ satisfy Garding type inequalities\footnote{\cite[Corollary~7.3]{AlvKordy1}.}.
Thus $\ker\bfDelta_1=\calH_2$, and the above decompositions restrict to \cinf~differential forms; i.e.,
\begin{align*}
{\calH}_1&=\ker\D_1\oplus\im d_1\oplus\im\delta_1\;,\\
\widetilde{\calH}_1&=\im\tilde d_1\oplus\im\tilde\delta_1\;,
\end{align*}
as topological vector spaces with the \cinf~topology, as well as with the restriction
of the $L^2$~norm topology.

\item[$($iv$)$] The space ${\calH}_2$ is of finite dimension, and the
inclusion ${\calH}_2\hookrightarrow{\calH}_1$ induces isomorphisms
$${\calH}_2^{u,v}\stackrel{\cong}{\lar} 
H^u\left({\calH}_1^{\cdot,v},d_1\right)
\cong H^u\left(\hat e_1^{\cdot,v}\right)
\cong H^u\left({\widehat E}_1^{\cdot,v}\right)\;.$$

\item[$($v$)$] We have $\tilde d_1^2=0$ and 
$H\left(\widetilde{\calH}_1,\tilde d_1\right)=0$.

\item[$($vi$)$] Each map $\widetilde{\calH}_1^{\cdot,v}\to\widetilde{\calH}_1^{\cdot,v}=\bar
o_1^{\cdot,v}\cong\bar0_1^{\cdot,v}$, defined by the canonical projection
$$
\overline{d_{0,1}\left(\Omega^{\cdot,v-1}\right)}
\oplus\overline{\delta_{0,-1}\left(\Omega^{\cdot,v}\right)}\lar
\overline{d_{0,1}\left(\Omega^{\cdot,v-1}\right)}/
d_{0,1}\left(\Omega^{\cdot,v-1}\right)\;,
$$
induces an isomorphism\footnote{The isomorphism
$H(\bar0_1)=0$ was originally shown  by X.~Masa \cite{Masa92}, as well as property~(vii), which
is a consequence.}
$$0=H^u\left(\widetilde{\calH}_1^{\cdot,v},d_1\right)
\cong H^u\left(\bar o_1^{\cdot,v}\right)\cong H^u\left(\bar0_1^{\cdot,v}\right)\;.$$

\item[$($vii$)$] All the following bigraded topological vector spaces are Hausdorff of finite
dimension and isomorphic to each other by maps that are either canonical or induced by the
projections $\pi_{u,v}$: $H({\hat e}_1)$, $H(e_1)$, $e_2$, $H({\widehat E}_1)$, $H(E_1)$ and
$E_2$.

\end{itemize} 
\end{Thm}

\lem{d10 P=P d10 P, ...}
The following properties are satisfied:
\begin{itemize}

\item[$($i$)$] We have
\begin{gather*}
\begin{alignat*}{3}
d_{1,0}P&=Pd_{1,0}P\;,&\qquad d_{1,0}\widetilde Q&=\widetilde Qd_{1,0}\widetilde Q\;,&\qquad
Qd_{1,0}&=Qd_{1,0}Q\;,\\
\widetilde Pd_{1,0}&=\widetilde Pd_{1,0}\widetilde P\;,&\qquad
\delta_{-1,0}Q&=Q\delta_{-1,0}Q\;,&\qquad
\delta_{-1,0}\widetilde P&=\widetilde P\delta_{-1,0}\widetilde P\;,
\end{alignat*}\\
P\delta_{-1,0}=P\delta_{-1,0}P\;,\qquad
\widetilde Q\delta_{-1,0}=\widetilde Q\delta_{-1,0}\widetilde Q\;.
\end{gather*}

\item[$($ii$)$] We have
$$
\widetilde Pd_{1,0}P=Qd_{1,0}\widetilde Q=\widetilde Q\delta_{-1,0}Q=
P\delta_{-1,0}\widetilde P=0\;.
$$

\end{itemize}
\elem

\prf
The equalities involving $d_{1,0}$ in property~(i) follow from~\refe{d0,12=...} since 
$$
P(\Omega)=\overline{d_{0,1}(\Omega)}\;,\quad\widetilde Q(\Omega)=\ker d_{0,1}\;.
$$
The other equalities in property~(i) are obtained by taking adjoints, and
property~(ii) is a direct consequence of property~(i).
\eprf

\lem{H1} The following properties are satisfied:
\begin{itemize}

\item[$($i$)$] The following operators on $\Omega$ define bounded operators on
$\bfOmega$:
\begin{alignat*}{4}
&\widetilde\Pi d_{1,0}\Pi\;,&\qquad&\Pi d_{1,0}\widetilde\Pi\;,&\qquad&
\widetilde\Pi\delta_{-1,0}\Pi\;,&\qquad&\Pi\delta_{-1,0}\widetilde\Pi\;,\\
&\widetilde Qd_{1,0}Q\;,&\qquad&Pd_{1,0}\widetilde P\;,&\qquad&  
\widetilde P\delta_{-1,0}P\;,&\qquad&
Q\delta_{-1,0}\widetilde Q\;.
\end{alignat*}

\item[$($ii$)$] The following operators on $\Omega$ define bounded operators on
$\bfOmega$ too:
\begin{alignat*}{4}
&\widetilde\Pi d\Pi\;,&\qquad&\Pi d\widetilde\Pi\;,&\qquad&
\widetilde\Pi_{\cdot,v+1}d\widetilde\Pi_{\cdot,v}\;,&\qquad& 
\widetilde\Pi_{\cdot,v-1}d\widetilde\Pi_{\cdot,v}\;,\\
&\widetilde\Pi\delta\Pi\;,&\qquad&\Pi\delta\widetilde\Pi\;,&\qquad&
\widetilde\Pi_{\cdot,v+1}\delta\widetilde\Pi_{\cdot,v}\;,&\qquad&
\widetilde\Pi_{\cdot,v-1}\delta\widetilde\Pi_{\cdot,v}\;.
\end{alignat*}

\item[$($iii$)$] We have 
\begin{alignat*}{2}
\dom\bfd_1&= L^2\calH_1\cap\dom\bfd\;,&\qquad
\dom\bfdelta_1&= L^2\calH_1\cap\dom\bfdelta\;,\\
\dom\tilde{\bfd}_1&= L^2\widetilde\calH_1\cap\dom\bfd\;,&\qquad
\dom\tilde{\bfdelta}_1&= L^2\widetilde\calH_1\cap\dom\bfdelta\;. 
\end{alignat*}

\end{itemize}
\elem

\prf Set $D_\perp=d_{1,0}+\delta_{-1,0}$.
Then, by Remark~3.7 and the proof of Lemma~7.2 in \cite{AlvKordy1}, the
operators 
$$[D_\perp,\Pi]\;,\quad
\widetilde\Pi_{\cdot,v-1}D_\perp\widetilde\Pi_{\cdot,v}\;,\quad
\left(\id-\widetilde\Pi_{\cdot,v}\right)D_\perp\widetilde\Pi_{\cdot,v}$$
on $\Omega$ define bounded operators on $\bfOmega$. This easily yields
property~(i). Now properties~(ii) and~(iii) follows 
from property~(i) since $d_{2,-1}$
and $\delta_{-2,1}$ are of order zero, and
$d_{0,1}$ and $\delta_{0,-1}$ vanish on ${\calH}_1$ and preserve each 
$\widetilde{\calH}_1^{\cdot,v}$.\eprf

\sec{L2 spectral seq}{$L^2$~spectral sequence}

\ssec{general L2}{General properties}

For a \cinf\ foliation \calF\ on a closed manifold $M$, what we call the
$L^2$~spectral sequence of \calF\ is
also a spectral sequence $(\bfE_k,\bfd_k)$ converging to the de~Rham
cohomology of $M$; in fact, it converges to the $L^2$~cohomology of $M$, but
both cohomologies are canonically isomorphic since $M$ is closed. Recall that $\bfd$ denotes the
closure of $d$ in $\bfOmega$. Also, let $\bfOmega_k$ be the
closure of $\Omega_k$ in $\bfOmega$, and consider the decreasing
filtration of the complex $(\dom\bfd,\bfd)$ by the differential subspaces
$\bfOmega_k\cap\dom\bfd$. We define $(\bfE_k,\bfd_k)$ to be the corresponding spectral
sequence. Since the inclusion $\Omega\hookrightarrow\dom\bfd$ obviously is a homomorphism of
filtered complexes, it induces a canonical homomorphism
$(E_k,d_k)\to(\bfE_k,\bfd_k)$ of spectral sequences. We point out that, by the compactness of $M$, 
the filtered complex $(\dom\bfd,\bfd)$ is well defined independently of any metric, and thus so is the
$L^2$~spectral sequence $(\bfE_k,\bfd_k)$.

Each $\bfE_1^{u,v}$ is a topological vector space with the topology
induced by the $L^2$~norm of $\bfOmega$, and consider the product
topology on $\bfE_1=\bigoplus_{u,v}\bfE_1^{u,v}$.

The notation 
$Z_k^{u,v}$ and $B_k^{u,v}$ of \refss{general diff} will be used for the spaces
involved in the definition of the differentiable spectral sequence of \calF, and
the corresponding spaces for the $L^2$~spectral sequence will be denoted by
$\bfZ_k^{u,v}$ and $\bfB_k^{u,v}$. We have 
\begin{align*}
\bfZ_k^{u,v}&=\bfOmega_u^{u+v}\cap
\bfd^{-1}\left(\bfOmega_{u+k}^{u+v+1}\right)\;,\\
\bfB_k^{u,v}&=\bfOmega_u^{u+v}\cap
\bfd\left(\bfOmega_{u-k}^{u+v-1}\cap\dom\bfd\right)\;,\\
\bfZ_\infty^{u,v}&=\bfOmega_u^{u+v}\cap \ker\bfd\;,\\    
\bfB_\infty^{u,v}&=\bfOmega_u^{u+v}\cap\im\bfd\;.
\end{align*}

As in the case of the differentiable spectral sequence, let $\pi_{u,v}:\bfOmega\to
\bfOmega^{u,v}$ be the canonical projection defined by the bigrading of $\bfOmega$; i.e.,
$\pi_{u,v}:\bfOmega\to\bfOmega^{u,v}$ is the continuous extension of
$\pi_{u,v}:\Omega\to\Omega^{u,v}$. Consider also the topological vector spaces
$$
\bfz_k^{u,v}=\pi_{u,v}\left(\bfZ_k^{u,v}\right)\;,\quad
\bfb_k^{u,v}=\pi_{u,v}\left(\bfB_k^{u,v}\right)\;,\quad
\bfe_k^{u,v}=\bfz_k^{u,v}/\bfb_{k-1}^{u,v}\;,\quad
\bfe_k=\bigoplus_{u,v}\bfe_k^{u,v}
$$
for $k=0,1,\ldots,\infty$, with the topology induced by the $L^2$~norm of $\bfOmega$. We
clearly have $\bfZ_k^{u,v}\cap\ker\pi_{u,v}=\bfZ_{k-1}^{u+1,v-1}$, and thus each projection
$\pi_{u,v}$ induces a continuous linear isomorphism
$\bfE_k^{u,v}\stackrel{\cong}{\to}\bfe_k^{u,v}$. Via these isomorphisms, the differential
$\bfd_k$ on $\bfE_k$ induces a differential on $\bfe_k$ that will be denoted by $\bfd_k$ as well. We
also have canonical continuous homomorphisms
$e_k^{u,v}\to\bfe_k^{u,v}$. 

In general, the $L^2$~spectral sequence is more difficult to deal with than the differentiable
spectral sequence. For example, we do not know whether the continuous linear isomorphism
$\bfE_1^{u,v}\stackrel{\cong}{\to}\bfe_1^{u,v}$, induced by $\pi_{u,v}$, is a homeomorphism with
this generality. Also, the useful expressions~\refe{Z0}--\refe{e1} do not hold for the
$L^2$~spectral sequence; indeed, for $r=u+v$, instead of~\refe{Z0}--\refe{Z1} we have
\begin{align} 
\bfZ_0^{u,v}&=\bfOmega_u^r\cap\dom\bfd\;,\label{e:bfZ0}\\
\bfB_0^{u,v}&=\bfd\left(\bfOmega_u^{r-1}\cap\dom\bfd\right)\;,\label{e:bfB0}\\
\bfZ_1^{u,v}&=\left(\left(\bfOmega^{u,v}\cap\ker\bfd_{0,1}\right)
+ \bfOmega_{u+1}^r\right)\cap\dom\bfd\;.\label{e:bfZ1}
\end{align}
Because of this reason, it will be useful to introduce the spaces
$$D^{u,v}=\pi_{u,v}\left(\bfOmega_u^r\cap\dom\bfd\right)\subset\bfOmega^{u,v}\;,\quad r=u+v\;,$$
which satisfy
\begin{equation}\label{e:D}
\left(V+ \bfOmega_{u+1}^r\right)\cap\dom\bfd=
\left(\left(V\cap D^{u,v}\right)+ \bfOmega_{u+1}^r\right)\cap\dom\bfd\;,\quad r=u+v\;.
\end{equation}
for any subspace $V\subset \bfOmega^{u,v}$.
%\begin{equation}\label{e:ker}
%\ker\left(\pi_{u,v}:\left(V+ \bfOmega_{u+1}^{u+v}\right)\cap\dom\bfd
%\lar D^{u,v}\right)=\bfOmega_{u+1}^{u+v}\cap\dom\bfd\;.
%\end{equation}

Observe that the canonical homomorphism $E_0^{u,v}\to \bfE_0^{u,v}$ is injective
with dense image because it is just the inclusion $Z_0^{u,v}\hookrightarrow
\bfZ_0^{u,v}$, whose image is dense by \refe{Z0} and \refe{bfZ0}. With this generality,
at least injectivity holds for $E_1\to \bfE_1$ too, as asserted by the following result.

\lem{E1 ar bfE1}
The canonical homomorphism $E_1\to \bfE_1$ is injective.
\elem

\prf For $r=u+v$ we have
\begin{align}
\lefteqn{\bfZ_0^{u+1,v-1}+ \bfB_0^{u,v}}\notag\\
&=\left(\bfOmega_{u+1}^r\cap\dom\bfd\right)
+\bfd\left(\bfOmega_u^{r-1}\cap\dom\bfd\right)\;,\quad\text{by \refe{bfB0} and
\refe{bfZ0}}\;,\notag\\
&=\left(\bfOmega_{u+1}^r
+\bfd\left(\bfOmega_u^{r-1}\cap\dom\bfd\right)\right)\cap\dom\bfd\;,\quad
\text{since}\quad\im\bfd\subset\dom\bfd\;,\notag\\  
&=\left(\bfd_{0,1}D^{u,v-1}+\bfOmega_{u+1}^r\right)\cap\dom\bfd\;.\label{e:bfZ0+bfB0}
\end{align} 
Then 
$$Z_1^{u,v}\cap\left(\bfZ_0^{u+1,v-1}+ \bfB_0^{u,v}\right)=
Z_0^{u+1,v-1}+B_0^{u,v}$$
by \refe{B0}, \refe{Z0} and \refe{Z1}, and the result follows.\eprf

\lem{D subset dom d01}
We have $D^{u,v}\subset\dom\bfd_{0,1}$.
\elem

\prf Take any $\alpha\in D^{u,v}$. For $r=u+v$, there exists some
$\beta\in\bfOmega_{u+1}^r$ such that $\alpha+\beta\in\dom\bfd$. So
$\pi_{u,v}\bfd(\alpha+\beta)$ is defined in $\bfOmega^{u,v}$. But 
$\pi_{u,v}\bfd(\alpha+\beta)=\bfd_{0,1}\alpha$ because $\alpha+\beta\in\bfOmega_u^r$.
\eprf 

\lem{bfe1}
We have
$$
\pi_{u,v}\left(\bfZ_1^{u,v}\right)=D^{u,v}\cap\ker\bfd_{0,1}\;,\quad
\pi_{u,v}\left(\bfB_0^{u,v}\right)=\bfd_{0,1}D^{u,v-1}\;,
$$
and thus
$$
\bfe_1^{u,v}=\frac{D^{u,v}\cap\ker\bfd_{0,1}}{\bfd_{0,1}D^{u,v-1}}\;.
$$
\elem

\prf For $r=u+v$, we have
\begin{alignat*}{2}
\pi_{u,v}\left(\bfZ_1^{u,v}\right)
&=\pi_{u,v}\left(\left(\left(\bfOmega^{u,v}\cap\ker\bfd_{0,1}\right)+
\bfOmega_{u+1}^r\right)\cap\dom\bfd\right)\;,&\qquad&\text{by \refe{bfZ1}}\;,\\ 
&=\pi_{u,v}\left(\left(\left(D^{u,v}\cap\ker\bfd_{0,1}\right)
+\bfOmega_{u+1}^r\right)\cap\dom\bfd\right)\;,&\qquad&\text{by \refe{D}}\;,\\
&=D^{u,v}\cap\ker\bfd_{0,1}\;,&&\\
\pi_{u,v}\left(\bfB_0^{u,v}\right)
&=\pi_{u,v}\left(\bfd\left(\bfOmega_u^{r-1}\cap\dom\bfd\right)\right)\;,&\qquad&
\text{by \refe{bfB0}}\;,\\ 
&=\pi_{u,v}\left(\left(\bfd_{0,1}D^{u,v-1}
+ \bfOmega_{u+1}^r\right)\cap\dom\bfd\right)\;,&\qquad&\text{by \refe{D}}\;,\\
&=\bfd_{0,1}D^{u,v-1}\;.\qed &&
\end{alignat*}
\renewcommand{\qed}{}
\eprf 

As for the differentiable spectral 
sequence, let $\bar\bfzero_1\subset\bfE_1$ and
$\bar\bfo_1\subset\bfe_1$ be the
closures of the corresponding trivial subspaces, which are bigraded subspaces with bigraded
quotients $\widehat\bfE_1=\bfE_1/\bar\bfzero_1$ and $\hat\bfe_1=\bfe_1/\bar\bfo_1$. 
\refl{bfe1} has the following direct consequence.

\cor{barbfo1,hatbfe1}
We have
\begin{gather*}
\bar\bfo_1^{u,v}=\frac{D^{u,v}\cap\cl_0\left(\bfd_{0,1}D^{u,v-1}\right)}
{\bfd_{0,1}D^{u,v-1}}
=\frac{D^{u,v}\cap\cl_0\left(d_{0,1}\Omega^{u,v-1}\right)}
{\bfd_{0,1}D^{u,v-1}}\;,\\
\hat\bfe_1^{u,v}=\frac{D^{u,v}\cap\ker\bfd_{0,1}}
{D^{u,v}\cap\cl_0\left(\bfd_{0,1}D^{u,v-1}\right)}
=\frac{D^{u,v}\cap\ker\bfd_{0,1}}{D^{u,v}
\cap\cl_0\left(d_{0,1}\Omega^{u,v-1}\right)}\;.
\end{gather*}
\ecor
 
The map $\bfd_1$, either on $\bfE_1$ or on $\bfe_1$, may not be continuous. So $\bar\bfzero_1$,
$\widehat\bfE_1$, $\bar\bfo_1$ and $\hat\bfe_1$ may not have canonical structures of bigraded 
complexes in general. However we shall show that this holds for Riemannian foliations in 
\refss{L2 spec seq Riem foln}.

\ssec{L2 spec seq Riem foln}{$L^2$~spectral sequence of Riemannian foliations}

\th{Ek=bfEk}
Let \calF\ be a Riemannian foliation on a closed manifold $M$.
Then the canonical map $E_k\to\bfE_k$ is injective with dense image for $k=0,1$, and is an
isomorphism of topological vector spaces for $k\geq 2$. In particular
$\bfE_k$ is Hausdorff of finite dimension for $k\geq2$.
\eth

The goal of this subsection is to prove \reft{Ek=bfEk}. Thus, from now on, assume \calF\ is a
Riemannian foliation. Since its statement is independent of any metric on $M$, we can take
a bundle-like metric on $M$ to prove it.

In \reft{Ek=bfEk}, the case $k=0$ is obvious, and the case $k=1$ follows directly from 
\refl{E1 ar bfE1} and the following lemma.

\lem{Z1 ar bfZ1}
The space $Z_1^{u,v}$ is dense in $\bfZ_1^{u,v}$.
\elem

\prf Since the orthogonal projection
$$\widetilde Q:\bfOmega^{u,v}\lar \bfOmega^{u,v}\cap\ker\bfd_{0,1}$$  
preserves smoothness on $M$, the result follows by \refe{Z1}
and \refe{bfZ1}.\eprf 

The proof of \reft{Ek=bfEk} for $k\geq2$ requires much more work than
\refl{Z1 ar bfZ1}. To establish this, we shall use the Hodge theoretic approach
to $e_1$ and $e_2$  from \refss{Hodge}, and a similar approach to $\bfe_1$ and $\bfe_2$. To begin with, we
show that $\bfd_1$ preserves $\bar\bfo_1$.

\lem{d(barbfo1) subset barbfo1}
We have $\bfd_1(\bar\bfo_1)\subset \bar\bfo_1$. 
\elem

\prf Take any $\alpha\in D^{u,v}\cap\cl_0\left(d_{0,1}\Omega^{u,v-1}\right)$,
and fix some $\beta\in \bfOmega_{u+1}^r$ with $\alpha+\beta\in\dom\bfd$, where $r=u+v$. We
know that $\pi_{u+1,v}\bfd(\alpha+\beta)\in D^{u+1,v}$ by \refl{bfe1}. On the other
hand, if $\bar\bfd_{0,1}$ and $\bar\bfd_{1,0}$ denote the extensions of
$d_{0,1}$ and $d_{1,0}$
to continuous maps
$\bfOmega\to W^{-1}\Omega$, we have
$$\pi_{u+1,v}\bfd(\alpha+\beta)=\bar\bfd_{1,0}\alpha+\bar\bfd_{0,1}\beta_1
\in\bar\bfd_{1,0}\alpha+\bar\bfd_{0,1}\bfOmega^{u+1,v-1}\;,$$
where $\beta_1=\pi_{u+1,v-1}\beta\in \bfOmega^{u+1,v-1}$, and
$$
\bar\bfd_{1,0}\alpha\in
\bar\bfd_{1,0}\left(\cl_0\left(d_{0,1}\Omega^{u,v-1}\right)\right)
\subset\cl_{-1}\left(d_{0,1}\Omega^{u+1,v-1}\right)\;.
$$
Hence
$$\pi_{u+1,v}\bfd(\alpha+\beta)\in
D^{u+1,v}\cap\cl_{-1}\left(d_{0,1}\Omega^{u+1,v-1}\right)
=D^{u+1,v}\cap\cl_0\left(d_{0,1}\Omega^{u+1,v-1}\right)$$
by \reft{leafwise Hodge}. Therefore the result follows by \refl{bfe1} and \refc{barbfo1,hatbfe1}.
\eprf

Now $\bar\bfo_1$ and $\hat\bfe_1$ canonically are bigraded complexes by
\refl{d(barbfo1) subset barbfo1}, and we have the short exact sequence
$$0\lar \bar\bfo_1
\lar \bfe_1
\lar \hat\bfe_1\lar0\;,$$
which induces long exact sequences
\begin{equation}\label{e:L2 long exact}
\cdots\lar H^u\left(\bar\bfo_1^{\cdot,v}\right)
\lar H^u\left(\bfe_1^{\cdot,v}\right)
\lar H^u\left(\hat\bfe_1^{\cdot,v}\right)
\lar H^{u+1}\left(\bar\bfo_1^{\cdot,v}\right)
\lar\cdots\;.
\end{equation}

\lem{D Hodge}
We have
$$
D^{u,v}\cap\ker\bfd_{0,1}=
\left(D^{u,v}\cap  L^2\calH_1\right)
\oplus\left(D^{u,v}\cap\cl_0\left(d_{0,1}\Omega^{u,v-1}\right)\right)
$$
as topological vector spaces, and moreover
$$
D^{u,v}\cap  L^2\calH_1= L^2\calH_1^{u,v}\cap\dom\bfd_1\;.
$$
\elem

\prf The inclusion ``$\supset$'' of the first equality
is obvious, and the inclusion ``$\supset$'' of the second equality follows from
\refl{H1}-(iii). 

To prove the inclusion ``$\subset$'' of the first equality,
by \refe{leafwise Hodge} it is enough
to prove that $\Pi\alpha\in D^{u,v}$ for all $\alpha\in
D^{u,v}\cap\ker\bfd_{0,1}$. This obviously holds if we prove
$\Pi\alpha\in\dom\bfd_1$ for every such an $\alpha$ since the inclusion
``$\supset$'' of the second equality is already proved. This also proves the
inclusion ``$\subset$'' of the second equality by taking $\alpha\in  L^2\calH_1$. 

Thus take any $\alpha\in D^{u,v}\cap\ker\bfd_{0,1}$. Then there is some
$\beta\in \bfOmega_{u+1}^r$ such that
$\alpha+\beta\in\dom\bfd$, where $r=u+v$. Write $\beta=\beta_1+\beta_2$ with $\beta_1\in
\bfOmega^{u+1,v-1}$ and $\beta_2\in \bfOmega_{u+2}^r$. Thus, since
$\alpha\in\ker\bfd_{0,1}$, we get 
$$\bfOmega^{u+1,v}\ni\Pi\pi_{u+1,v}\bfd(\alpha+\beta)
=\Pi\left(\bar\bfd_{1,0}\alpha+\bar\bfd_{0,1}\beta_1\right)
=\Pi\bar\bfd_{1,0}\alpha\;.$$
Here we consider 
%$\bar\bfd_{1,0}\alpha+\bar\bfd_{0,1}\beta_1$ as an element of $W^{-1}\Omega$,
%and 
$\Pi$ and $\pi_{u+1,v}$ as bounded operators on $W^{-1}\Omega$. But
$$\Pi\bar\bfd_{1,0}\alpha
=\Pi\bar\bfd_{1,0}\Pi\alpha+\Pi\bar\bfd_{1,0}P\alpha$$
because $Q\alpha=0$, and 
$$\Pi\bar\bfd_{1,0}P\alpha=\Pi\widetilde{P}\bar\bfd_{1,0}P\alpha\in \bfOmega$$
by \refl{H1}-(i). Therefore $\Pi\bar\bfd_{1,0}\Pi\alpha\in\bfOmega$, yielding
$\Pi\alpha\in\dom\bfd_1$ as desired.\eprf

\cor{hatbfe1 Hodge}
The inclusions $ L^2\calH_1^{u,v}\cap\dom\bfd_1\hookrightarrow D^{u,v}\cap\ker\bfd_{0,1}$ induce an
isomorphism of $(\dom\bfd_1,\bfd_1)$ onto the quotient complex $\hat\bfe_1$, which is also an
isomorphism of topological vector spaces.
\ecor

\prf This follows from \refc{barbfo1,hatbfe1} and \refl{D Hodge}.\eprf 

\cor{H(hatbfe1) Hodge}
Each inclusion $\calH_2^{u,v}\hookrightarrow D^{u,v}\cap\ker\bfd_{0,1}$ induces an isomorphism
$H^u\left(\hat\bfe_1^{\cdot,v}\right)\cong\calH_2^{u,v}$
of topological vector spaces. In particular $H(\hat\bfe_1)$ is Hausdorff of finite dimension.
\ecor

\prf This follows from
\refc{hatbfe1 Hodge} and
\reft{H1}-(iii),(iv).\eprf 

The canonical homomorphism $e_1\to \bfe_1$ is obviously continuous. Hence
it induces homomorphisms of complexes $\bar o_1\to\bar\bfo_1$ and 
$\hat e_1\to\hat\bfe_1$, and homomorphisms 
$H(\bar o_1)\to H(\bar\bfo_1)$ and  
$H(\hat e_1)\to H(\hat\bfe_1)$ in cohomology.

\cor{H(hate1)=H(hatbfe1)}
The canonical map  $H(\hat e_1)\to H(\hat\bfe_1)$ is an isomorphism of topological vector spaces.
\ecor

\prf This follows from \reft{H1}-(iv) and \refc{H(hatbfe1) Hodge}.\eprf 

We also need a Hodge theoretic study of certain complex whose cohomology is
isomorphic to
$H\left(\bar\bfo_1\right)$. To simplify notation let\footnote{This notation is 
used in \cite{AlvKordy1} for the \cinf~versions of these complexes.}
$$
\calZ_v=\bigoplus_u\bfZ_1^{u,v}\;,\quad
\calB_v=\bigoplus_u\left(\bfZ_0^{u-1,v+1}+\bfB_0^{u,v}\right)\;,
$$
which are subcomplexes of $(\dom\bfd,\bfd)$.
Then
\begin{equation}\label{e:bf01}
\bar{\bf0}_1^{\cdot,v}=\cl_0(\calB_v)/\calB_v\;.
\end{equation}
Observe that $\calZ_{v-1}\subset\calB_v$.

\lem{Sergiescu}
The quotient complex ${\calB}_v/{\calZ}_{v-1}$ is acyclic.
Thus the quotient map
$\cl_0(\calB_v)/\calZ_{v-1}\to\cl_0({\calB}_v)/{\calB}_v=\bar{\bf0}_1^{\cdot,v}$
induces an  isomorphism in cohomology.
\elem

\prf The result follows from \refe{bfZ1} and \refe{bfZ0+bfB0} with easy arguments
(see Lemma~2.5 in \cite{Sergiescu87}  and Lemma~7.4 in \cite{AlvKordy1}). \eprf

Set 
\begin{align*}
\widetilde{\bfOmega}^{u,v}&=\bfOmega^{u,v}+ \bfOmega^{u+1,v-1}\;,\\
\tilde\pi_{u,v}&=\pi_{u,v}+\pi_{u+1,v-1}:\bfOmega\lar\widetilde{\bfOmega}^{u,v}\;,\\
\widetilde{D}^{u,v}&=\tilde\pi_{u,v}\left(\bfOmega_u^r\cap\dom\bfd\right)\;,\quad r=u+v\;.
\end{align*}
We have
$$
\cl_0(\calB_v)\cap\ker\tilde\pi_{u,v}\subset\calZ_{v-1}\cap\ker\tilde\pi_{u,v}\;.
$$
Hence, for each topological vector space
$$
\tilde{\bfe}_1^{u,v}
=\frac{\tilde\pi_{u,v}\left(\cl_0(\calB_v)\right)}{\tilde\pi_{u,v}\left(\calZ_{v-1}\right)}\;,
$$
the projection $\tilde\pi_{u,v}$ induces a continuous linear isomorphism
\begin{equation}\label{e:cl0(B)/Z=tildebfe1}
\cl_0(\calB_v^r)/\calZ_{v-1}^r\stackrel{\cong}{\lar}
\tilde{\bfe}_1^{u,v}\;,\quad r=u+v\;.
\end{equation} 
Let $\tilde\bfd_1$ be the operator on $\tilde{\bfe}_1=\bigoplus_{u,v}\tilde{\bfe}_1^{u,v}$
that corresponds to the differential operator on the quotient complex
$\cl_0(\calB_v)/\calZ_{v-1}$ by the above isomorphisms. Observe that $\tilde\bfd_1$ is given as
follows: if
$\alpha\in\cl_0(\calB_v)$, and $[\tilde\pi_{u,v}\alpha]\in\tilde{\bfe}_1^{u,v}$ denotes the
class defined by $\tilde\pi_{u,v}\alpha$, then
$\tilde\bfd_1[\tilde\pi_{u,v}\alpha]=[\tilde\pi_{u+1,v}\bfd\alpha]$.

The spaces $\widetilde{D}^{u,v}$ and $D^{u,v}$ have similar properties. For instance, for any
subspace $V\subset\widetilde{\bfOmega}^{u,v}$ we have
\begin{equation}\label{e:tilde D}
\left(V+ \bfOmega_{u+2}^r\right)\cap\dom\bfd=
\left(\left(V\cap\widetilde{D}^{u,v}\right)
+ \bfOmega_{u+2}^r\right)\cap\dom\bfd\;,\quad r=u+v\;.
\end{equation}
%\begin{equation}\label{e:ker tilde}
%\ker\left(\tilde\pi_{u,v}:\left(V+ \bfOmega_{u+2}^{u+v}\right)\cap\dom\bfd
%\lar\widetilde{D}^{u,v}\right)=\bfOmega_{u+2}^{u+v}\cap\dom\bfd\;.
%\end{equation}

\lem{tildebfe1}
For $r=u+v$, we have
\begin{gather*}
\tilde\pi_{u,v}\left(\cl_0(\calB_v)\right)=\widetilde{D}^{u,v}\cap\cl_0({\calB}_v)\;,\\
\tilde\pi_{u,v}\left(\calZ_{v-1}\right)=D^{u+1,v-1}\cap{\calZ}_{v-1}
=D^{u+1,v-1}\cap\ker\bfd_{0,1}\;,
\end{gather*}
and thus
$$
\tilde{\bfe}_1^{u,v}
=\frac{\widetilde{D}^{u,v}\cap\cl_0({\calB}_v)}{D^{u+1,v-1}\cap\ker\bfd_{0,1}}\;.
$$ 
\elem

\prf This easily follows from \refe{D} and \refe{tilde D}.\eprf 

The following result and  \refl{D Hodge} are similar, as well as their proofs.

\lem{tilde D Hodge}
We have
$$\widetilde{D}^{u,v}\cap\cl_0({\calB}_v)=
\left(\widetilde{D}^{u,v}\cap L^2\widetilde\calH_1\right)
\oplus\left(D^{u+1,v-1}\cap\ker\bfd_{0,1}\right)$$
as topological vector spaces, and moreover
$$\widetilde{D}^{u,v}\cap L^2\widetilde\calH_1= L^2\widetilde\calH_1^{u,v}
\cap\dom\tilde\bfd_1\;.$$
\elem

\prf The inclusion ``$\supset$'' of the first equality is
obvious, and the inclusion ``$\supset$'' of the second equality follows from
\refl{H1}-(iii). 

To prove the inclusion ``$\subset$'' of the first equality,
by \refe{leafwise Hodge} it is enough
to prove that $\widetilde\Pi\alpha\in\widetilde{D}^{u,v}$ for all $\alpha\in
\widetilde{D}^{u,v}\cap\cl_0({\calB}_v)$. This obviously holds if we prove
$\widetilde\Pi\alpha\in\dom\tilde\bfd_1$ for every such an $\alpha$ since the
inclusion ``$\supset$'' of the second equality is already proved. This also
proves the inclusion ``$\subset$'' of the second equality by taking $\alpha\in
 L^2\widetilde\calH_1$. 

Thus take any $\alpha\in \widetilde{D}^{u,v}\cap\cl_0({\calB}_v)$. Then there
is some
$\beta\in \bfOmega_{u+2}^r$ such that
$\alpha+\beta\in\dom\bfd$, where $r=u+v$. Write $\alpha=\alpha_1+\alpha_2$ with $\alpha_1\in
\bfOmega^{u,v}$ and $\alpha_2\in \bfOmega^{u+1,v-1}$. So, since
$\alpha\in\cl_0({\calB}_v)$ and
$\bar\bfd\beta\in\ker\widetilde\Pi_{\cdot,v}$, where $\widetilde\Pi_{\cdot,v}$ is
considered as a projection in $W^{-1}\Omega$, we get
$$\widetilde{\bfOmega}^{u+1,v}\ni
\widetilde\Pi_{\cdot,v}\bfd(\alpha+\beta)
=\widetilde\Pi_{\cdot,v}\bar\bfd\alpha\;.$$ 
%Here we consider $\bar\bfd\alpha$ as an element of $W^{-1}\Omega$, and
%$\widetilde\Pi_{\cdot,v}$ as bounded operator on $W^{-1}\Omega$.  
But
$$\widetilde\Pi_{\cdot,v}\bar\bfd\alpha=
\widetilde\Pi_{\cdot,v}\bar\bfd\widetilde\Pi_{\cdot,v}\alpha
+\widetilde\Pi_{\cdot,v}\bar\bfd\Pi\alpha_2
+\widetilde\Pi_{\cdot,v}\bar\bfd\widetilde\Pi_{\cdot,v-1}\alpha_2$$ because
$\alpha\in
\widetilde{D}^{u,v}\cap\cl_0({\calB}_v)$, and
$$\widetilde\Pi_{\cdot,v}\bar\bfd\Pi\alpha_2
+\widetilde\Pi_{\cdot,v}\bar\bfd\widetilde\Pi_{\cdot,v-1}\alpha_2\in \bfOmega$$
by \refl{H1}-(ii). Therefore
$\widetilde\Pi_{\cdot,v}\bar\bfd\widetilde\Pi_{\cdot,v}\alpha\in\bfOmega$, yielding
$\widetilde\Pi_{\cdot,v}\alpha\in\dom\tilde\bfd_1$ as desired.\eprf 

Consider the projection 
$$
\widetilde{D}^{u,v}\cap\cl_0({\calB}_v)\lar
\widetilde{D}^{u,v}\cap L^2\widetilde\calH_1= L^2\widetilde\calH_1^{u,v}\cap\dom\tilde\bfd_1
$$ 
defined by \refl{tilde D Hodge}, which is obviously an orthogonal projection.

\cor{tildebfe1 Hodge}
The inclusions
$ L^2\widetilde\calH_1^{u,v}\cap\dom\tilde\bfd_1\hookrightarrow
\widetilde{D}^{u,v}\cap\cl_0({\calB}_v)$
induce an isomorphism
$\left(\dom\tilde{\bfd}_1,\tilde{\bfd}_1\right)\stackrel{\cong}{\to}
\left(\tilde{\bfe}_1,\tilde{\bfd}_1\right)$
of bigraded complexes and topological vector spaces.
\ecor

\prf This follows from \refls{tildebfe1}{tilde D Hodge}.\eprf 

\cor{H(barbfo1)=0}
We have $H\left(\bar\bfo_1\right)=0$.
\ecor

\prf This follows from \refe{bf01}, \refe{cl0(B)/Z=tildebfe1}, \refl{Sergiescu},
\refc{tildebfe1 Hodge} and \reft{H1}-(v). \eprf 

\cor{H(bfe_1)=H(hatbfe_1)}
The canonical map $H(\bfe_1)\to H(\hat\bfe_1)$ is an
isomorphism of topological vector spaces. In particular $H(\bfe_1)$ is Hausdorff of finite
dimension.
\ecor

\prf The canonical map $H(\bfe_1)\to H(\hat\bfe_1)$ is a linear isomorphism by
\refc{H(barbfo1)=0} and the exactness of \refe{L2 long exact}. Moreover it is obviously
continuous. Then it is also an homeomorphism because
$H(\hat\bfe_1)$ is a Hausdorff topological vector space of finite dimension. \eprf 

\cor{H(e_1)=H(bfe_1)}
The canonical map $H(e_1)\to H(\bfe_1)$ is an isomorphism of topological vector spaces.
\ecor

\prf By the commutativity of the diagram
$$
\begin{CD}
H(e_1) & @>>> &H(\bfe_1)\\
@VVV & & @VVV\\
H(\hat e_1) &@>>> &H(\hat\bfe_1)\;,
\end{CD}
$$
where all maps are canonical, the result follows directly from \reft{H1}-(vii),
\refcs{H(hate1)=H(hatbfe1)}{H(bfe_1)=H(hatbfe_1)}.
\eprf

\cor{E2=bfE2)}
The canonical map $E_2\to\bfE_2$ is an isomorphism of topological vector spaces.
\ecor

\prf Consider the compositions 
$$E_2\lar H(E_1)\lar H(e_1)\;,\quad\bfE_2\lar H(\bfE_1)\lar H(\bfe_1)\;,$$ 
where the first map of each composition is canonical, and the second one is 
canonically induced by the projections $\pi_{u,v}$. The first composition is an
isomorphism of topological vector spaces by \reft{H1}-(vii), and we know that the second
composition is a continuous linear isomorphism (\refss{general L2}). Then the second composition is also an homeomorphism because
$H(\bfe_1)$ is Hausdorff of finite dimension by \refc{H(bfe_1)=H(hatbfe_1)}. So the result
follows from
\refc{H(e_1)=H(bfe_1)} and the commutativity of the diagram
$$
\begin{CD}
E_2 & @>>> &H(e_1)\\
@VVV & & @VVV\\
\bfE_2 &@>>> &H(\bfe_1)\;,
\end{CD}
$$
where the horizontal arrows denote the above compositions, and the vertical arrows denote
canonical maps. 
\eprf

Now \reft{Ek=bfEk} for $k\geq2$ follows from \refc{E2=bfE2)} because the canonical map
$(E_k,d_k)\to(\bfE_k,\bfd_k)$ is a homomorphism of spectral sequences.

\sec{L2 spec seq small eigenvalues}{$L^2$~spectral sequence and small eigenvalues}

\ssec{main results}{Main results}

Let \calF\ be a \cinf~foliation on a closed manifold $M$ with a Riemannian metric $g$, and consider the
family of metrics $g_h$, $h>0$, which were defined in \refe{gh} and give rise to the adiabatic
limit. As in \refs{intro}, let $\D_{g_h}$ denote the Laplacian on $\Omega$ defined by
$g_h$, and 
$$0\leq\lambda_0^r(h)\leq\lambda_1^r(h)\leq\lambda_2^r(h)\leq\cdots$$
its spectrum on $\Omega^r$, taking multiplicity into
account. The following result suggests that, with this generality, 
the number of small eigenvalues of $\D_h$ may be more related with the
$L^2$~spectral sequence than with the differentiable one. Nevertheless, so far we do not know
about the relevance its hypothesis for non-Riemannian foliations.

\th{bfZk-1+bfZinfty closed in bfZk}
Let \calF\ be a \cinf~foliation on a closed Riemannian manifold. If 
$\bfZ_{k-1}^{u+1,v-1}+\bfZ_\infty^{u,v}$ is closed in $\bfZ_k^{u,v}$ for all
$u,v$, with $r=u+v$, then
$$\dim \bfE_k^r\leq
\sharp\,\left\{i\ \left|\ \lambda_i^r(h)\in O\left(h^{2k}\right.\right)
\quad\text{as}\quad h\downarrow 0\right\}
$$
for all $r$.
\eth

The following more understandable result is a direct
consequence of \reft{bfZk-1+bfZinfty closed in bfZk} because
$$\frac{\bfZ_\ell^{u,v}}{\bfZ_{l-1}^{u+1,v-1}+\bfZ_\infty^{u,v}}$$
is a quotient of $\bfE_\ell^{u,v}$.

\cor{bfEkr Hausdorff of finite dim}
Let \calF\ be a \cinf~foliation on a closed Riemannian manifold. If $\bfE_k$ is
Hausdorff of finite dimension, then
$$\dim\bfE_\ell^r\leq
\sharp\,\left\{i\ \left|\ \lambda_i^r(h)\in O\left(h^{2l}\right.\right)
\quad\text{as}\quad h\downarrow 0\right\}\;,\quad \ell\geq k\;.
$$
\ecor

\rem{leq}
Observe that, by \reft{Ek=bfEk}, \refc{bfEkr Hausdorff of finite dim} holds for Riemannian foliations and
$k=2$, and inequality ``$\leq$'' of \refe{small eigenvalues k} in \refmt{small eigenvalues} follows.
\erem

The proof of \reft{bfZk-1+bfZinfty closed in bfZk} is given in \refss{proof}, and its two main ingredients
are described in Sections~\ref{ss:spectral distribution function} and~\ref{ss:decomposition}: the variational
formula of the spectral distribution function used by Gromov-Shubin, and the direct sum decomposition for
general spectral sequences.

\ssec{spectral distribution function}{Spectral distribution function}

For a closed Riemannian manifold $(M,g)$, let $N^r(\lambda)$ denote the spectral
distribution function of the Laplacian $\D$ on $\Omega^r$; i.e. $N^r(\lambda)$
is the number of eigenvalues of $\D$ on $\Omega^r$ which are $\leq \lambda$, 
taking multiplicity into account. 
Recall that $\bfOmega$ denotes the Hilbert space of square integrable
differential forms with the inner product induced by $g$, and $\bfd$ the closure of the de~Rham derivative
$d$ in $\bfOmega$. Let $\bar\bfd:\dom\bfd/\ker\bfd\to\bfOmega$
denote the map induced by $\bfd$, and consider the quotient Hilbert norm on
$\bfOmega/\ker\bfd$. The following variational expression of $N^r(\lambda)$ is a consequence of the Hodge
decomposition of $\bfOmega$.

\begin{Prop}[{Gromov-Shubin \cite{GromovShubin}}]\label{p:GromovShubin}
We have $$N^r(\lambda)=F^{r-1}(\lambda)+\beta^r+F^r(\lambda)\;,$$ 
where $\beta^r$ is the $r$th Betti number of $M$, and
$$F^r(\lambda)=\sup_L\dim L\;,$$
with $L$ ranging over the closed subspaces of $\dom\bfd/\ker\bfd$ satisfying
$$\norm{\bar\bfd\zeta}\leq\sqrt{\lambda}\,\norm{\zeta}\quad\text{for
all}\quad\zeta\in L\;.$$
\end{Prop}

Now take again a \cinf~foliation \calF\ on $M$. Then, for each metric $g_h$
of the family~\refe{gh} that gives rise to the adiabatic limit, the spectral distribution
function of $\D_{g_h}$ will be denoted by $N_h^r(\lambda)$, and decomposes as
$$N_h^r(\lambda)=F_h^{r-1}(\lambda)+\beta^r+F_h^r(\lambda)\;,$$
according to \refp{GromovShubin}. 

Suppose \calF\ is of codimension $q$, and let $\norm{\ }_h$ be the norm induced
by $g_h$ on $\bfOmega$. The following equality will be also used to prove \reft{bfZk-1+bfZinfty closed in
bfZk}:
\begin{equation}\label{e:| |h}
\norm{\omega}_h=h^{-q/2}h^u\,\norm{\omega}
\quad\text{if}\quad\omega\in\Omega^{u,v}\;.
\end{equation}
This follows from two observations. First, if the
metrics induced by $g$ and $g_h$ on
$\bigwedge TM^\ast$ are also denoted by $g$ and $g_h$, then
$g_h=h^{2u}g$ on forms with transverse degree $u$. And second,  assuming $M$ is oriented, the volume forms
$\mu$ and $\mu_h$, induced by $g$ and $g_h$, satisfy 
$\mu_h=h^{-q}\mu$ since volume forms are of transverse degree $q$.

By using \refp{GromovShubin} in the same spirit of \cite{GromovShubin}, we could prove that
the asymptotics of the $\lambda_i^r(h)$, as $h\downarrow0$, are \cinf~homotopy invariants of \calF\
(with respect to the appropriate definition of homotopy between foliations). However, for our
purposes in this paper, it will be enough to prove that the asymptotics of the
$\lambda_i^r(h)$ are independent of the choice of the given metric $g$ on $M$. This will not be used to prove
\reft{bfZk-1+bfZinfty closed in bfZk} but will play an important role to finish the proof of \refmt{small
eigenvalues} in \refss{rescaling}. Such independence of $g$ is proved in the following way. Let $g'$ be
another metric on $M$ with corresponding 1-parameter family of metrics $g'_h$, and
let $\norm{\ }'$ and $\norm{\ }'_h$ denote the corresponding norms on
$\bfOmega$. Compactness of
$M$ implies the existence of some $C>0$ such that 
$$
C^{-1}\norm{\omega}\leq\norm{\omega}'\leq C\norm{\omega}
$$ 
for all $\omega\in\bfOmega$, yielding
\begin{equation}\label{e:| |'h}
C^{-1}\norm{\omega}_h\leq\norm{\omega}'_h\leq C\norm{\omega}_h
\end{equation}
for all $\omega\in\bfOmega$ and $h>0$ by \refe{| |h}. Let $N_h^{'r}(\lambda)$ be the spectral distribution
function of
$\D_{g'_h}$ on $\Omega^r$, and let
$$
N_h^{\prime r}(\lambda)=
F_h^{\prime r-1}(\lambda)+\beta^r+F_h^{\prime r}(\lambda)
$$ 
be its decomposition according to \refp{GromovShubin}. Then
$$
F_h^{\prime r}(C^{-4}\lambda)\leq F_h^r(\lambda)\leq 
F_h^{\prime r}(C^4\lambda)
$$
for all $\lambda\geq0$ $h>0$ by \refe{| |'h} and the definition of $F_h^r$ and $F_h^{\prime r}$.
Thus
\begin{equation}\label{e:N'rh}
N_h^{\prime r}(C^{-4}\lambda)\leq F_h^r(\lambda)\leq 
N_h^{\prime r}(C^4\lambda)\;,
\end{equation}
yielding the metric independence of the asymptotics of the $\lambda_i^r(h)$.

\ssec{decomposition}{Direct sum decomposition of spectral sequences} 

In this subsection we consider the general setting where $(E_k,d_k)$ is the spectral sequence induced by an
arbitrary complex $(\calA,d)$ with a finite  decreasing filtration
$$\calA=\calA_0\supset\calA_1\supset\cdots
\supset\calA_q\supset\calA_{q+1}=0$$
by differential subspaces. 

\lem{decomposition of calAr}
The following properties are satisfied:
\begin{itemize}

\item[$($i$)$] There is a $($non-canonical$)$ isomorphism
\begin{align*}
\calA^r&\cong E_\infty^r\oplus\bigoplus_\ell\left(\left(E_\ell^r\cap\im d_\ell\right)
\oplus\frac{E_\ell^r}{E_\ell^r\cap\ker d_\ell}\right)\\
&=\bigoplus_{u+v=r}\left(E_\infty^{u,v}\oplus\bigoplus_\ell\left(\left(E_\ell^{u,v}\cap\im
d_\ell\right)\oplus\frac{E_\ell^{u,v}}{E_\ell^{u,v}\cap\ker d_\ell}\right)\right)\;.
\end{align*}

\item[$($ii$)$] The isomorphism in $($i$)$ can be chosen so that $\calA_k^r$ corresponds to
$$
\bigoplus_{u\geq k,\ u+v=r}
\left(E_\infty^{u,v}\oplus\bigoplus_\ell\left(\left(E_\ell^{u,v}\cap\im
d_\ell\right)\oplus\frac{E_\ell^{u,v}}{E_\ell^{u,v}\cap\ker d_\ell}\right)\right)\;.
$$

\item[$($iii$)$] The isomorphism in $($i$)$ can be
chosen so that the only possibly non-trivial components of the operator corresponding to $d$ by
$($i$)$ are the isomorphisms
$$\bar d_\ell:\frac{E_\ell^{u,v}}{E_\ell^{u,v}\cap\ker d_\ell}\lar
E_\ell^{u+\ell,v-\ell+1}\cap\im d_\ell$$
canonically defined by $d_\ell$.

\end{itemize}
\elem

Before proving \refl{decomposition of calAr}, we state three corollaries that will be needed in the proof
of \refp{GromovShubin}.

\cor{decomposition of Ekr}
There is a $($non-canonical$)$ isomorphism
$$
E_k^r\cong
E_\infty^r\oplus\bigoplus_{l\geq k}\left(\left(E_\ell^r\cap\im d_\ell\right)\oplus
\frac{E_\ell^r}{E_\ell^r\cap\ker d_\ell}\right)\;.
$$
\ecor

\prf
This is a direct consequence of \refl{decomposition of calAr}.
\eprf

Let
$$
m_k^r=\dim\bigoplus_{l\geq k}\frac{E_\ell^r}{E_\ell^r\cap\ker d_\ell}\;.
$$

\cor{m}
We have 
$$\dim E_k^r=m_k^{r-1}+H^r(\calA,d)+m_k^r\;.$$
\ecor

\prf This follows from \refc{decomposition of Ekr} since each
$d_\ell$ induces isomorphisms
$$\frac{E_\ell^r}{E_\ell^r\cap\ker d_\ell}\cong E_\ell^{r+1}\cap\im d_\ell\;. \qed$$
\renewcommand{\qed}{}\eprf 

\cor{L}
For $r=u+v$, there is a subspace $L_k^{u,v}\subset\calA^r/(\calA^r\cap\ker d)$ such that:
\begin{itemize}

\item[$($i$)$] We have
$$
\frac{Z_k^{u,v}+(\calA^r\cap\ker d)}{\calA^r\cap\ker d}
=L_k^{u,v}\oplus\frac{Z_{k-1}^{u+1,v-1}+(\calA^r\cap\ker d)}{\calA^r\cap\ker d}
$$
as vector spaces. In particular $\bar d\left(L_k^{u,v}\right)\subset
\calA_{u+k}^{r+1}$.

\item[$($ii$)$] The direct sum $L_k^r=\bigoplus_{u+v=r}L_k^{u,v}$ makes sense in
$\calA^r/(\calA^r\cap\ker d)$, and we have $\dim L_k^r=m_k^r$.

\end{itemize}
\ecor

\prf
>From \refl{decomposition of calAr} we get a (non-canonical) isomorphism
\begin{equation}\label{e:A/ker d}
\frac{\calA^r}{\calA^r\cap\ker d}
\cong\bigoplus_\ell\frac{E_\ell^r}{E_\ell^r\cap\ker d_\ell}\;,
\end{equation}
Then let $L_k^{u,v}$ be the subspace of $\calA^r/(\calA^r\cap\ker d)$ that corresponds to
$$
\bigoplus_{l\geq k}\frac{E_\ell^{u,v}}{E_\ell^{u,v}\cap\ker d_\ell}
$$
by \refe{A/ker d}. Then property~(i) easily follows from \refl{decomposition of calAr}, and
property~(ii) is obvious; in fact, $L_k^r$ corresponds to
$$
\bigoplus_{l\geq k}\frac{E_\ell^r}{E_\ell^r\cap\ker d_\ell}
$$
by \refe{A/ker d}.
\eprf

\rem{L}
By \refc{L}-(i), the canonical isomorphism
\begin{equation}\label{e:Z/Z=Z/ker d}
\frac{Z_k^{u,v}}{Z_{k-1}^{u+1,v-1}+Z_\infty^{u,v}}\stackrel{\cong}{\lar}
\frac{Z_k^{u,v}+(\calA^r\cap\ker d)}
{Z_{k-1}^{u+1,v-1}+(\calA^r\cap\ker d)}
\end{equation}
yields
$$
\frac{Z_k^{u,v}}{Z_{k-1}^{u+1,v-1}+Z_\infty^{u,v}}\cong L_k^{u,v}\;.
$$
\erem

When applying \refc{L} to the $L^2$~spectral sequence of a \cinf~foliation, the subspaces
$L_k^r\subset\dom\bfd/\ker\bfd$ of \refc{L} will be the spaces $L$ needed to apply \refp{GromovShubin}.

The rest of this section will be devoted to prove \refl{decomposition of calAr}. To begin
with, we have \cite{McClearly}
\begin{align*}
E_\ell^{u,v}\cap d_\ell(E_\ell)
&=\frac{B_\ell^{u,v}+Z_{\ell-1}^{u+1,v-1}}
{Z_{\ell-1}^{u+1,v-1}+B_{\ell-1}^{u,v}}\;,\\
E_\ell^{u,v}\cap\ker d_\ell
&=\frac{Z_{\ell+1}^{u,v}+Z_{\ell-1}^{u+1,v-1}}
{Z_{\ell-1}^{u+1,v-1}+B_{\ell-1}^{u,v}}\;.
\end{align*}
So
\begin{align}
E_\ell^{u,v}\cap d_\ell(E_\ell)&\cong
\frac{B_\ell^{u,v}}{B_{\ell+1}^{u+1,v-1}+B_{\ell-1}^{u,v}}\;,
\label{e:dl(El)=B/B}\\
\frac{E_\ell^{u,v}}{E_\ell^{u,v}\cap\ker d_\ell}&\cong
\frac{Z_\ell^{u,v}}{Z_{\ell-1}^{u+1,v-1}+Z_{\ell+1}^{u,v}}
\label{e:El/ker dl=Z/Z}
\end{align}
canonically. Here, isomorphism~\refe{El/ker dl=Z/Z} is obvious, and 
\refe{dl(El)=B/B} follows since 
$$
B_{\ell-1}^{u,v}\subset B_\ell^{u,v}\;,\quad 
B_\ell^{u,v}\cap Z_{\ell-1}^{u+1,v-1}=B_{\ell+1}^{u+1,v-1}\;.
$$

Consider the following chain of inclusions for
$0\leq u\leq q$ and $r=u+v$:
\begin{equation}\label{e:inclusions}
\left.\begin{array}{c}
\calA_{u+1}^r\subset\calA_{u+1}^r+B_0^{u,v}
\subset\calA_{u+1}^r+B_1^{u,v}\subset\cdots\\[6pt]
\cdots\subset\calA_{u+1}^r+B_\infty^{u,v}
\subset\calA_{u+1}^r+Z_\infty^{u,v}\subset\cdots\\[6pt]
\cdots\subset\calA_{u+1}^r+Z_2^{u,v}\subset
\calA_{u+1}^r+Z_1^{u,v}\subset\calA_u^r\,.
\end{array}\right\}
\end{equation}
The inclusions in \refe{inclusions} 
have the following quotients:
\begin{align}
\frac{\calA_{u+1}^r+B_\ell^{u,v}}{\calA_{u+1}^r+B_{\ell-1}^{u,v}}
&\cong\frac{B_\ell^{u,v}}{B_{\ell+1}^{u+1,v-1}+B_{\ell-1}^{u,v}}\;,
\label{e:B/B}\\
\frac{\calA_{u+1}^r+Z_\infty^{u,v}}{\calA_{u+1}^r+B_\infty^{u,v}}
&\cong\frac{Z_\infty^{u,v}}{Z_\infty^{u+1,v-1}+B_\infty^{u,v}}
=E_\infty^{u,v}\;,
\label{e:Z/B}\\
\frac{\calA_{u+1}^r+Z_\ell^{u,v}}{\calA_{u+1}^r+Z_{\ell+1}^{u,v}}
&\cong\frac{Z_\ell^{u,v}}{Z_{\ell-1}^{u+1,v-1}+Z_{\ell+1}^{u,v}}\;,
\label{e:Z/Z}
\end{align}
where these isomorphisms are canonical because
$$
\begin{array}{ll}
B_{\ell-1}^{u,v}\subset B_\ell^{u,v}\;,&\quad 
B_\ell^{u,v}\cap\calA_{u+1}^r=B_{\ell+1}^{u+1,v-1}\;,\\[6pt]
B_\infty^{u,v}\subset Z_\infty^{u,v}\;,&\quad 
Z_\infty^{u,v}\cap\calA_{u+1}^r=Z_\infty^{u+1,v-1}\;,\\[6pt]
Z_{\ell+1}^{u,v}\subset Z_\ell^{u,v}\;,&\quad 
Z_\ell^{u,v}\cap\calA_{u+1}^r=Z_{\ell-1}^{u+1,v-1}\;. 
\end{array}
$$

The direct sum decomposition in property~(i) will depend on
the  choice of linear complements for the inclusions in \refe{inclusions}:
\begin{align*}
\calA_{u+1}^r+B_\ell^{u,v}&=U_\ell^{u,v}\oplus
\left(\calA_{u+1}^r+B_{\ell-1}^{u,v}\right)\;,\\
\calA_{u+1}^r+Z_\infty^{u,v}&=V^{u,v}\oplus
\left(\calA_{u+1}^r+B_\infty^{u,v}\right)\;,\\
\calA_{u+1}^r+Z_\ell^{u,v}&=W_\ell^{u,v}\oplus
\left(\calA_{u+1}^r+Z_{\ell+1}^{u,v}\right)\;.
\end{align*}
On the one hand, since the chains in \refe{inclusions} form a filtration of $\calA^r$ when varying
$u$, we have
\begin{equation}\label{e:A=bigoplus UVW}
\calA^r=\bigoplus_{u+v=r}\left(V^{u,v}\oplus\bigoplus_\ell\left(U_\ell^{u,v}\oplus
W_\ell^{u,v}\right)\right)
\end{equation}
as vector space. On the other hand, according to the canonical isomorphisms \refe{B/B},~\refe{Z/B}
and~\refe{Z/Z}, the spaces $U_\ell^{u,v}$, $V^{u,v}$ and $W_\ell^{u,v}$ can be chosen so that 
\begin{equation}\label{e:UVW}
U_\ell^{u,v}\subset B_\ell^{u,v}\;,\quad V^{u,v}\subset Z_\infty^{u,v}\;,\quad 
W_\ell^{u,v}\subset Z_\ell^{u,v}\;,
\end{equation}
yielding direct sum decompositions
\begin{align}
B_\ell^{u,v}&=U_\ell^{u,v}\oplus
\left(B_{\ell+1}^{u+1,v-1}+B_{\ell-1}^{u,v}\right)\;,\label{e:U}\\
Z_\infty^{u,v}&=V^{u,v}\oplus
\left(Z_\infty^{u+1,v-1}+B_\infty^{u,v}\right)\;,\label{e:V}\\
Z_\ell^{u,v}&=W_\ell^{u,v}\oplus
\left(Z_{\ell-1}^{u+1,v-1}+Z_{\ell+1}^{u,v}\right)\;.\label{e:W}
\end{align}
Hence
\begin{align}
U_\ell^{u,v}&\cong E_\ell^{u,v}\cap\im d_\ell\;,
\label{e:U=im dl}\\
V^{u,v}&\cong E_\infty^{u,v}\;,\label{e:V=Einfty}\\
W_\ell^{u,v}&\cong\frac{E_\ell^{u,v}}{E_\ell^{u,v}\cap\ker d_\ell}
\label{e:W=el/ker dl}
\end{align}
by \refe{dl(El)=B/B},~\refe{El/ker dl=Z/Z} and~\refe{U}--\refe{W}. 
Therefore property~(i) follows from~\refe{A=bigoplus UVW} and~\refe{U=im dl}--\refe{W=el/ker dl}.

Property~(ii) follows from \refe{A=bigoplus UVW} because
$$U_\ell^{u,v}\;,\;V^{u,v}\;,\;W_\ell^{u,v}\subset\calA_u^r\;.$$

Now property~(iii) is obviously equivalent to the existence of $U_\ell^{u,v}$, $V^{u,v}$ and
$W_\ell^{u,v}$ as above satisfying
\begin{equation}\label{e:dUV=0, dW subset U}
d\left(U_\ell^{u,v}\right)=d\left(V^{u,v}\right)=0\;,\quad
d\left(W_\ell^{u,v}\right)=U_\ell^{u+\ell,v-\ell+1}\;.
\end{equation}
The first equality of \refe{dUV=0, dW subset U} holds by \refe{UVW}. We shall also check that,
once the $W_\ell^{u,v}$ is given satisfying \refe{W}, the $U_\ell^{u,v}$ defined by 
\refe{dUV=0, dW subset U} satisfies \refe{U}. This follows because
$d$ canonically induces a map
$$
\hat d_\ell:
\frac{Z_\ell^{u,v}}{Z_{\ell-1}^{u+1,v-1}+Z_{\ell+1}^{u,v}}\lar
\frac{B_\ell^{u+\ell,v-\ell+1}}{B_{\ell+1}^{u+\ell+1,v-\ell}+B_{\ell-1}^{u+\ell,v-\ell+1}}\;,
$$
which corresponds to the isomorphism $\bar d_\ell$
via  \refe{El/ker dl=Z/Z} and \refe{dl(El)=B/B}.
So $\hat d_\ell$ is an isomorphism as well, and thus the above $U_\ell^{u,v}$
satisfies \refe{U} as desired. This finishes the proof of \refl{decomposition of calAr}.

\ssec{proof}{Proof of \reft{bfZk-1+bfZinfty closed in bfZk}}

Assume $\bfZ_{k-1}^{u+1,v-1}+\bfZ_\infty^{u,v}$ is closed in $\bfZ_k^{u,v}$ for all
$u,v$. We shall need the following abstract result.

\lem{Hilbert}
Let $L$ be a real complete metrizable topological vector space, and $V,W\subset L$ linear subspaces. If
$V\cap W=0$, $V$ is closed in $L$, and $W$ is closed in $V+W$, then $V+W=V\oplus W$ as topological
vector spaces.
\elem

\prf We have $(V+W)\cap\overline{W}=W$ since $W$ is closed in $V+W$, yielding
$V\cap\overline{W}=0$ because
$V\cap W=0$. So 
$V+\overline{W}=V\oplus\overline{W}$ as topological vector spaces because all spaces
involved are closed subspaces of $L$ (see for instance \cite[Corollary~3 of Theorem~2.1, Chapter~III,
page~78]{Schaefer}). Now the result follows easily.
\eprf

\lem{ker d+Omega closed in Z+ker d+Omega}
For $u+v=r$, the space 
$(\bfOmega^r\cap\ker\bfd)+\bfOmega_{u+1}^r$ is a closed subspace of
$\bfZ_k^{u,v}+(\bfOmega^r\cap\ker\bfd)+\bfOmega_{u+1}^r$.
\elem

\prf The space $\bfOmega^r\cap\ker\bfd$ is closed in $\bfOmega$ since $\bfd$ is a closed operator,
and thus so is its subspace $\bfZ_\infty^{u,v}=\bfOmega_u\cap(\bfOmega^r\cap\ker\bfd)$. Hence
$\bfOmega^r\cap\ker\bfd=V\oplus\bfZ_\infty^{u,v}$
as Hilbert spaces, where $V$ is the orthogonal complement of $\bfZ_\infty^{u,v}$ in 
$\bfOmega^r\cap\ker\bfd$; in particular $V$ is closed in $\bfOmega$ too. Obviously,
\begin{align*}
(\bfOmega^r\cap\ker\bfd)+\bfOmega_{u+1}^r&=V+\bfZ_\infty^{u,v}+\bfOmega_{u+1}^r\;,\\
\bfZ_k^{u,v}+(\bfOmega^r\cap\ker\bfd)+\bfOmega_{u+1}^r
&=V+\bfZ_k^{u,v}+\bfOmega_{u+1}^r\;.
\end{align*}
On the other hand we clearly have
\begin{align*}
\bfZ_\infty^{u,v}+\bfOmega_{u+1}^r
&=\bfOmega_u\cap\left((\bfOmega^r\cap\ker\bfd)+\bfOmega_{u+1}^r\right)\;,\\
\bfZ_k^{u,v}+\bfOmega_{u+1}^r
&=\bfOmega_u\cap\left(\bfZ_k^{u,v}+(\bfOmega^r\cap\ker\bfd)+\bfOmega_{u+1}^r\right)\;,
\end{align*}
and thus $\bfZ_\infty^{u,v}+\bfOmega_{u+1}^r$ and $\bfZ_k^{u,v}+\bfOmega_{u+1}^r$ are
respectively closed in $(\bfOmega^r\cap\ker\bfd)+\bfOmega_{u+1}^r$ and 
$\bfZ_k^{u,v}+(\bfOmega^r\cap\ker\bfd)+\bfOmega_{u+1}^r$. Therefore \refl{Hilbert} yields
\begin{align*}
(\bfOmega^r\cap\ker\bfd)+\bfOmega_{u+1}^r
&=V\oplus\left(\bfZ_\infty^{u,v}+\bfOmega_{u+1}^r\right)\;,\\
\bfZ_k^{u,v}+(\bfOmega^r\cap\ker\bfd)+\bfOmega_{u+1}^r
&=V\oplus\left(\bfZ_k^{u,v}+\bfOmega_{u+1}^r\right)\;,
\end{align*}
as topological vector spaces, and the result follows. 
\eprf

\rem{technical difficulty}
In the proof of \refl{ker d+Omega closed in Z+ker d+Omega}, the existence of $V$ so that
$\bfOmega^r\cap\ker\bfd=V\oplus\bfZ_\infty^{u,v}$ as Hilbert spaces is the technical difficulty we were not
able to solve without using square integrable differential forms; that is, we do not know if 
$\Omega^r\cap\ker d=V\oplus Z_\infty^{u,v}$ as topological vector spaces for some subspace $V$. This is the
whole reason of introducing the $L^2$~spectral sequence in this paper. 

Also, observe that the formula of Gromov-Shubin uses square integrable
differential forms. Thus it can be more easily related to the 
$L^2$~spectral sequence than to the differentiable one. Though this is a minor problem that could be easily
solved in the setting of
\cinf~differential forms.
\erem

We shall use the notation
$$
X_u^r=\bigoplus_{a\leq u}\Omega^{a,r-a}\;,\quad\rho_u^r=\sum_{a\leq u}\pi_{a,r-a}:\bfOmega^r\lar X_u^r\;.
$$
With respect to the inner product in $\bfOmega$ induced by $g$ or any $g_h$, the space $X_u^r$ is the
orthogonal complement of $\bfOmega_{u+1}^r$ in $\bfOmega^r$, and $\rho_u^r$ is an orthogonal projection.

\cor{rho(ker d) closed in rho(Z+ker d)}
For $u+v=r$, the space 
$\rho_u^r(\bfOmega^r\cap\ker\bfd)$ is closed in
$\rho_u^r\left(\bfZ_k^{u,v}+(\bfOmega^r\cap\ker\bfd)\right)$.
\ecor

\prf This follows from \refl{ker d+Omega closed in Z+ker d+Omega} since we clearly have
\begin{align*}
(\bfOmega^r\cap\ker\bfd)+\bfOmega_{u+1}^r
&=\rho_u^r(\bfOmega^r\cap\ker\bfd)\oplus\bfOmega_{u+1}^r\;,\\
\bfZ_k^{u,v}+(\bfOmega^r\cap\ker\bfd)+\bfOmega_{u+1}^r
&=\rho_u^r\left(\bfZ_k^{u,v}+(\bfOmega^r\cap\ker\bfd)\right)\oplus\bfOmega_{u+1}^r\;,
\end{align*}
as topological vector spaces. \eprf

Recall that $\bar\bfd:\dom\bfd/\ker\bfd\to\im\bfd$ denotes the map induced by $\bfd$, and
let $L_k^{u,v}$ and $L_k^r$ be the spaces introduced in
\refc{L} in \refss{decomposition} for the particular case of the $L^2$~spectral sequence of \calF.

\lem{d zeta}
We have
$$\norm{\bar\bfd\zeta}_h\leq h^{-q/2}h^{u+k}\norm{\bar\bfd\zeta}$$
for all $\zeta\in L_k^{u,v}$ and $0<h\leq 1$.
\elem

\prf This follows directly from \refc{L} and \refe{| |h}.\eprf

Let $\norm{\cdot}$ and $\norm{\cdot}_h$ also stand for the quotient
Hilbert norms on $\bfOmega/\ker\bfd$ induced by the norms $\norm{\cdot}$ and $\norm{\cdot}_h$ on $\bfOmega$,
respectively. In particular we have the restrictions of $\norm{\cdot}$ and $\norm{\cdot}_h$ to each subspace
$L_k^{u,v}\subset\bfOmega/\ker\bfd$.

\lem{C'K}
For each subspace $K\subset L_k^{u,v}$ of finite dimension 
there is some $C'_K>0$,
depending on $K$, such that
$$h^{-q/2}h^u\norm{\zeta}\leq C'_K\norm{\zeta}_h$$
for all $\zeta\in K$ and $0<h\leq 1$.
\elem

\prf Let $u+v=r$. The restriction $\rho_u^r:\bfZ_k^{u,v}+(\bfOmega^r\cap\ker\bfd)\to X_u^r$ induces a
homomorphism 
$$
\bar\rho_u^r:\frac{\bfZ_k^{u,r-u}+(\bfOmega^r\cap\ker\bfd)}{\bfOmega^r\cap\ker\bfd}
\lar\frac{X_u^r}{\rho_u^r\left(\bfOmega^r\cap\ker\bfd\right)}\;.
$$
We clearly have
\begin{equation}\label{e:ker rho}
\ker\bar\rho_u^r=\frac{\bfOmega_{u+1}^r+\left(\bfOmega^r\cap\ker\bfd\right)}
{\bfOmega^r\cap\ker\bfd}\;.
\end{equation}
So $\bar\rho_u^r$ induces a continuous linear isomorphism
$$
\frac{\bfZ_k^{u,v}+(\bfOmega^r\cap\ker\bfd)}
{\bfZ_{k-1}^{u+1,v-1}+(\bfOmega^r\cap\ker\bfd)}
\stackrel{\cong}{\lar}
\im\bar\rho_u^r=\frac{\rho_u^r\left(\bfZ_k^{u,v}+(\bfOmega^r\cap\ker\bfd)\right)}
{\rho_u^r(\bfOmega^r\cap\ker\bfd)}\;.
$$
Observe that $\im\bar\rho_u^r$ is a Hausdorff topological vector space by 
\refc{rho(ker d) closed in rho(Z+ker d)}, and thus $\norm{\cdot}$ and $\norm{\cdot}_h$ induce norms on
$\im\bar\rho_u^r$ that will be also denoted by $\norm{\cdot}$ and $\norm{\cdot}_h$, respectively. By
\refe{ker rho} and \refc{L}, the homomorphism
$\bar\rho_u^r$ restricts to an injection
$\bar\rho_u^r:L_k^{u,v}\to\im\bar\rho_u^r$.
Since $\rho_u^r$ is an orthogonal projection for any metric $g_h$, we easily get
\begin{equation}\label{e:rho}
\norm{\bar\rho_u^r\zeta}_h\leq\norm{\zeta}_h\quad\text{for all}\quad\zeta\in L_k^{u,v}\;.
\end{equation}
Here, we use the norm on $\im\bar\rho_u^r$ in the left hand side of \refe{rho}, and the
norm on $\bfOmega/\ker\bfd$ in its right hand side. Observe that, by \refe{| |h},
$$
h^{-q/2}h^u\norm{\omega}\leq\norm{\omega}_h\quad\text{for all}\quad
\omega\in X_u^r\quad\text{and}\quad 0<h\leq 1\;,
$$
yielding
\begin{equation}\label{e:im rho}
h^{-q/2}h^u\norm{\xi}\leq\norm{\xi}_h\quad\text{for all}\quad
\xi\in\im\bar\rho_u^r\quad\text{and}\quad 0<h\leq 1\;.
\end{equation}
Moreover, since $K$ is of finite dimension, $\im\bar\rho_u^r$ is Hausdorff, and the
restriction $\bar\rho_u^r:L_k^{u,v}\to\im\bar\rho_u^r$ is injective, we get the existence of
some $C'_K>0$ so that
\begin{equation}\label{e:C'K}
\norm{\zeta}\leq C'_K\norm{\bar\rho_u^r\zeta}\quad\text{for all}\quad\zeta\in K\;.
\end{equation}
So
\begin{alignat*}{2}
h^{-q/2}h^u\norm{\zeta}&\leq C'_Kh^{-q/2}h^u\norm{\bar\rho_u^r\zeta}\;,&&
\qquad\text{by~\refe{C'K}}\;,\\
&\leq C'_K\norm{\bar\rho_u^r\zeta}_h\;,&&\qquad\text{by~\refe{im rho}}\;,\\
&\leq C'_K\norm{\zeta}_h\;,&&\qquad\text{by~\refe{rho}}\;,
\end{alignat*}
for all $\zeta\in K$ and $0<h\leq 1$ as desired.
\eprf 

\cor{d zeta}
For each subspace $K\subset L_k^r$ of finite dimension 
there is some $C_K>0$, depending on $K$, such that
$$\norm{\bar\bfd\zeta}_h\leq C_Kh^k\norm{\zeta}_h$$
for all $\zeta\in K$ and $0<h\leq 1$.
\ecor

\prf Since $K$ is of finite dimension,
there is some constant $C''_K$, depending on $K$, so that 
\begin{equation}\label{e:C''K}
\norm{\bar\bfd\zeta}\leq C''_K\norm{\zeta}\qquad\text{for all}\quad\zeta\in
K\;.
\end{equation}

Because $L_k^r=\bigoplus_{u+v=r}L_k^{u,v}$, any finite
dimensional subspace $K\subset L_k^r$ is contained in the sum of finite dimensional
subspaces
$K^{u,v}\subset L_k^{u,v}$, $u+v=r$. Therefore we can assume $K$ is contained in some 
$L_k^{u,v}$ with $u+v=r$. Then, for
$\zeta\in K$ and $0<h\leq 1$, we have
\begin{alignat*}{2}
\norm{\bar\bfd\zeta}_h&\leq h^{-q/2}h^{u+k}\norm{\bar\bfd\zeta}\;,&&
\qquad\text{by \refl{d zeta}}\;,\\
&\leq C''_Kh^{-q/2}h^{u+k}\norm{\zeta}\;,&&\qquad\text{by \refe{C''K}}\;,\\
&\leq C'_KC''_Kh^k\norm{\zeta}_h\;,&&\qquad\text{by \refl{C'K}}\;,
\end{alignat*}
and the result follows with $C_K=C'_KC''_K$.\eprf 

Now the proof of \reft{bfZk-1+bfZinfty closed in bfZk} can be finished as follows. If
$m_k^r<\infty$, then \refc{d zeta} holds for $K=L_k^r$, and thus
$$F_h^r\left(C_\ell h^{2k}\right)\geq m_k^r\;.$$
Therefore, in this case, \reft{bfZk-1+bfZinfty closed in bfZk} follows from
\refc{m} and \refp{GromovShubin}. 

If $m_k^r=\infty$, choose any sequence of finite dimensional subspaces
$K_i\subset L_k^r$ so that $\dim K_i\uparrow\infty$. Then \refc{d zeta}
gives a sequence $C_i>0$ such that 
$$F_h^r\left(C_ih^{2k}\right)\geq\dim K_i$$
for $0<h\leq1$. Hence \reft{bfZk-1+bfZinfty closed in bfZk} also follows in this case by
\refc{m} and \refp{GromovShubin}.

\sec{asymptotics}{Asymptotics of eigenforms}

In the whole of this section, \calF\ is assumed to be a Riemannian foliation and the metric bundle-like.

\ssec{Hodge nested seq}{The Hodge theoretic nested sequence}

So far we have constructed bigraded subspaces $\calH_1,\calH_2\subset\Omega$,
which are respectively isomorphic to $\hat e_1,e_2$ as bigraded topological vector spaces by
\reft{leafwise Hodge} and \reft{H1}-(iv). We continue constructing subspaces  
$\calH_k\subset\Omega$ and isomorphisms $e_k\cong\calH_k$ by induction on $k$ as follows. Suppose
we have constructed $\calH_k$ and an explicit isomorphism $e_k\cong\calH_k$ for some $k\geq2$.
Then the homomorphism $d_k$ corresponds to some homomorphism on $\calH_k$ that
will be denoted by $d_k$ as well. Thus $\calH_k$ becomes a finite
dimensional complex. Let $\delta_k$ be the adjoint of $d_k$ on the finite dimensional Hilbert
space $\calH_k$, and set $\D_k=d_k\delta_k+\delta_kd_k$ and $\calH_{k+1}=\ker\D_k=\ker
d_k\cap\ker\delta_k$. We have the orthogonal decomposition
$$\calH_k=\calH_{k+1}\oplus\im d_k\oplus\im\delta_k\;,$$
yielding
$$e_{k+1}\cong H(e_k,d_k)\cong H(\calH_k,d_k)\cong\calH_{k+1}\;,$$
which completes the induction step. So $(\calH_k,d_k)$ is, by definition, some kind of a Hodge
theoretic version of the sequence $(\hat e_1,d_1),(e_2,d_2),(e_3,d_3),\ldots$, 
and thus of the sequence $(\widehat E_1,d_1),(E_2,d_2),(E_3,d_3),\ldots$ 
as well by \reft{leafwise Hodge} and \reft{H1}-(vii). Furthermore each
$\D_k$ is bihomogeneous of bidegree $(0,0)$, and thus $\calH_k$ inherits the bigrading from
$\bfOmega$, which clearly corresponds to the bigrading of
$E_k$ and $e_k$. Observe that the nested sequence
$$\Omega\supset\calH_1\supset\calH_2\supset\calH_3
\supset\calH_4\supset\cdots$$
stabilizes at most at the $(q+1)$th step since so does $E_k$. Then its
final term $\calH_{q+1}=\calH_{q+2}=\cdots$ will be denoted by $\calH_\infty$, and
we have $E_\infty\cong e_\infty\cong\calH_\infty$.

We shall need a better understanding of the new terms $\calH_k$ for $k>2$. Precisely, we shall
use the following result.

\th{calHk}
Let $k\geq3$ and $\omega\in\calH_2^{u,v}$. Then $\omega\in\calH_k^{u,v}$ if and only if there are
sequences $\alpha_i=\sum_{a>0}\alpha_i^a$ and $\beta_i=\sum_{a>0}\beta_i^a$,
where $\alpha_i^a\in\Omega^{u+a,v-a}$ and $\beta_i^a\in\Omega^{u-a,v+a}$, such that
$$\pi_{u+a,v-a+1}d(\omega+\alpha_i)\lar0\;,\quad\pi_{u-a,v+a-1}\delta(\omega+\beta_i)\lar0$$ 
strongly in $\bfOmega$ for $0<a<k$.
\eth

The rest of this section will be devoted to prove \reft{calHk}. 
To begin with, the nested sequence $\calH_k$ is most properly a Hodge theoretic version of
another sequence of bigraded topological complexes $(\hat e_{1,k},d_k)$, which are defined as follows by
induction on $k\ge1$. First, let $\hat e_{1,1}=\hat e_1$ and $\hat e_{1,2}=H(\hat e_1)$ with the
induced topology in cohomology. We have an explicit isomorphism $e_2\cong\hat e_{1,2}$ of
bigraded topological vector spaces given by
\reft{H1}-(vii). Now suppose that, for some fixed
$k\ge2$, we have defined $\hat e_{1,k}$ with an explicit isomorphism
$e_k\cong\hat e_{1,k}$ of bigraded topological vector spaces. Then $\hat e_{1,k}$ becomes a
topological complex via this isomorphism, and define $\hat e_{1,k+1}=H(\hat e_{1,k})$.
Furthermore the composition $e_{k+1}\cong H(e_k)\cong\hat e_{1,k+1}$ is an explicit isomorphism
of bigraded topological vector spaces. 

\lem{hat e1k}
For $k\ge1$, we have a canonical isomorphism
\begin{equation}\label{e:hate1,ku,v=...}
\hat e_{1,k}^{u,v}\cong
\frac{z_k^{u,v}+\overline{b_0^{u,v}}}{b_{k-1}^{u,v}+\overline{b_0^{u,v}}}
\end{equation}
of topological vector spaces. Moreover, for $k\ge2$, the above isomorphism $e_k^{u,v}\cong\hat
e_{1,k}^{u,v}$ corresponds to the canonical map
\begin{equation}\label{e:zku,v/bk-1u,v to ...}
\frac{z_k^{u,v}}{b_{k-1}^{u,v}}
\lar\frac{z_k^{u,v}+\overline{b_0^{u,v}}}{b_{k-1}^{u,v}+\overline{b_0^{u,v}}}
\end{equation}
when applying \refe{hate1,ku,v=...}.
\elem

\prf 
The result is proved by induction on $k$.
First, the case $k=1$ is trivial. 

Second, the kernel and the image of $d_1$ in $e_1^{u,v}$
respectively are
$z_2^{u,v}/b_0^{u,v}$ and
$b_1^{u,v}/b_0^{u,v}$, whose canonical projections in $\hat
e_1^{u,v}=z_1^{u,v}/\overline{b_0^{u,v}}$ are
\begin{equation}\label{e:ker and im of d1}
\frac{z_2^{u,v}+\overline{b_0^{u,v}}}{\overline{b_0^{u,v}}}\;,\quad
\frac{b_1^{u,v}+\overline{b_0^{u,v}}}{\overline{b_0^{u,v}}}\;,
\end{equation}
yielding the canonical isomorphism \refe{hate1,ku,v=...} for $k=2$. Since the
isomorphism $e_2^{u,v}\cong\hat e_{1,2}^{u,v}$ is canonically defined, it corresponds to the
canonical map \refe{zku,v/bk-1u,v to ...} for $k=2$.

Now assume the result holds for $k=\ell\geq2$ and we prove it for $k=\ell+1$. The kernel and the image
of $d_\ell$ in $e_\ell^{u,v}$ respectively are $z_{\ell+1}^{u,v}/b_{\ell-1}^{u,v}$ and
$b_\ell^{u,v}/b_{\ell-1}^{u,v}$, whose images by the canonical isomorphism
\refe{zku,v/bk-1u,v to ...} for $k=\ell$ are
\begin{equation}\label{e:ker and im of dl}
\frac{z_{\ell+1}^{u,v}+\overline{b_0^{u,v}}}{b_{\ell-1}^{u,v}+\overline{b_0^{u,v}}}\;,\quad
\frac{b_\ell^{u,v}+\overline{b_0^{u,v}}}{b_{\ell-1}^{u,v}+\overline{b_0^{u,v}}}\;.
\end{equation}
These spaces respectively correspond to the kernel and the image of $d_\ell$ in $\hat e_{1,\ell}^{u,v}$ by
\refe{hate1,ku,v=...}, yielding  the canonical isomorphism \refe{hate1,ku,v=...} for $k=\ell+1$.
Again, because the isomorphism $e_{\ell+1}^{u,v}\cong\hat e_{1,\ell+1}^{u,v}$ is canonically defined, it is
given by the canonical map \refe{zku,v/bk-1u,v to ...} for $k=\ell+1$.
\eprf

We shall consider each isomorphism \refe{hate1,ku,v=...} as an equality from now on.

For $k\geq1$, let $\Pi_k$ denote the orthogonal projections $\bfOmega\to\calH_k$; in particular,
$\Pi_1=\Pi$ with this notation. Let also $P_0=P$, $Q_0=Q$ and, for
$k\geq1$, let $P_k$ and $Q_k$ be the orthogonal projections of $\bfOmega$ onto
$d_k(\calH_k)$ and $\delta_k(\calH_k)$. Finally let $\bar P_k=\sum_{0\leq\ell\leq k}P_\ell$ and
$\bar Q_k=\sum_{0\leq\ell\leq k}Q_\ell$ for $k\geq0$. 

\lem{calHk}
For $k\ge1$, $\Pi_k$ induces an isomorphism $\hat
e_{1,k}^{u,v}\stackrel{\cong}{\to}\calH_k^{u,v}$, whose composition with the canonical
isomorphism $e_k^{u,v}\stackrel{\cong}{\to}\hat e_{1,k}^{u,v}$ is the above isomorphism
$e_k^{u,v}\cong\calH_k^{u,v}$.  
\elem

\prf 
Observe that the first part of the statement means that we have an orthogonal decomposition
\begin{equation}\label{e:zk+...=calHk oplus bk-1+...}
z_k^{u,v}+\overline{b_0^{u,v}}=
\calH_k^{u,v}\oplus\left(b_{k-1}^{u,v}+\overline{b_0^{u,v}}\right)\;.
\end{equation}

Again the result follows by induction on $k$. We have an orthogonal decomposition
\begin{equation}\label{e:z1=calH1 oplus ...}
z_1^{u,v}=\calH_1^{u,v}\oplus\overline{b_0^{u,v}}
\end{equation}
by \reft{leafwise Hodge}. Thus the isomorphism $\hat e_{1,1}=\hat
e_1^{u,v}\stackrel{\cong}{\to}\calH_1^{u,v}$ is induced by the orthogonal projection $\Pi_1$ onto
$\calH_1$. On the other hand, the kernel and image of
$d_1$ in $\calH_1^{u,v}$ respectively correspond by this isomorphism to the kernel and image of $d_1$ on
$\hat e_1^{u,v}$, which are respectively given by \refe{ker and im of d1}. So the kernel and image of
$d_1$ in $\calH_1^{u,v}$ are the orthogonal projections
$\Pi_1\left(z_2^{u,v}+\overline{b_0^{u,v}}\right)$ and
$\Pi_1\left(b_1^{u,v}+\overline{b_0^{u,v}}\right)$, respectively. Hence, by definition,
$\calH_2^{u,v}$ is the orthogonal complement of
$\Pi_1\left(b_1^{u,v}+\overline{b_0^{u,v}}\right)$ in
$\Pi_1\left(z_2^{u,v}+\overline{b_0^{u,v}}\right)$, which
is equal to the orthogonal complement of $b_1^{u,v}+\overline{b_0^{u,v}}$ in
$z_2^{u,v}+\overline{b_0^{u,v}}$ by \refe{z1=calH1 oplus ...} since
$$\overline{b_0^{u,v}}\subset b_1^{u,v}+\overline{b_0^{u,v}}
\subset z_2^{u,v}+\overline{b_0^{u,v}}\subset z_1^{u,v}\;.$$
Thus the result follows for $k=2$.

Now suppose the statement holds for $k=\ell\geq2$. Then, via the isomorphism $\hat
e_{1,\ell}^{u,v}\stackrel{\cong}{\to}\calH_\ell^{u,v}$ induced by $\Pi_\ell$, the kernel and image of
$d_\ell$ in $\calH_\ell^{u,v}$ respectively correspond to the kernel and image of $d_\ell$ in $\hat
e_{1,\ell}^{u,v}$, which are given in \refe{ker and im of dl}. So the kernel and image of
$d_\ell$ in $\calH_\ell^{u,v}$ are the orthogonal projections
$\Pi_\ell\left(z_{\ell+1}^{u,v}+\overline{b_0^{u,v}}\right)$ and
$\Pi_\ell\left(b_\ell^{u,v}+\overline{b_0^{u,v}}\right)$, respectively. Hence, by definition,
$\calH_{\ell+1}^{u,v}$ is the orthogonal complement of
$\Pi_\ell\left(b_\ell^{u,v}+\overline{b_0^{u,v}}\right)$ in
$\Pi_\ell\left(z_{\ell+1}^{u,v}+\overline{b_0^{u,v}}\right)$, which
is equal to the orthogonal complement of $b_\ell^{u,v}+\overline{b_0^{u,v}}$ in
$z_{\ell+1}^{u,v}+\overline{b_0^{u,v}}$ by \refe{zk+...=calHk oplus bk-1+...} for $k=\ell$ since
$$b_{\ell-1}^{u,v}+\overline{b_0^{u,v}}\subset b_\ell^{u,v}+\overline{b_0^{u,v}}
\subset z_{\ell+1}^{u,v}+\overline{b_0^{u,v}}\subset z_\ell^{u,v}+\overline{b_0^{u,v}}\;.$$
Thus the result follows for $k=\ell+1$.
\eprf

\rem{summarize} 
The inverse of the isomorphism
$\hat e_{1,k}^{u,v}\stackrel{\cong}{\to}\calH_k^{u,v}$ in \refl{calHk} is obviously induced by the
inclusion $\calH_k^{u,v}\hookrightarrow z_k^{u,v}+\overline{b_0^{u,v}}$. So we can summarize
\refls{hat e1k}{calHk} by saying that, for $k\ge2$, the isomorphism
$e_k^{u,v}\cong\calH_k^{u,v}$ is given by the diagram
\begin{equation}\label{e:summarize}
e_k^{u,v}=\frac{z_k^{u,v}}{b_{k-1}^{u,v}}\stackrel{\cong}{\lar}
\hat e_{1,k}^{u,v}=\frac{z_k^{u,v}+\overline{b_0^{u,v}}}{b_{k-1}^{u,v}+\overline{b_0^{u,v}}}
\stackrel{\cong}{\longleftarrow}\calH_k^{u,v}\;,
\end{equation}
where both isomorphisms are canonically induced by inclusions.
\erem

\rem{hate1,k}
In general, we have $z_k^{u,v}\neq\calH_k^{u,v}\oplus b_{k-1}^{u,v}$ because
$\calH_k^{u,v}\not\subset z_k^{u,v}$, but the orthogonal decomposition~\refe{zk+...=calHk oplus
bk-1+...} always holds. This is the reason the nested sequence $\calH_k$ is a Hodge theoretic
version of the sequence $(\hat e_{1,k},d_k)$ better than of the sequence $(\hat
e_1,d_1),(e_2,d_2),(e_3,d_3),\ldots$.
\erem

The following proposition is the key result
to prove \reft{calHk}.

\prop{dk in calHk}
Let $\omega\in\calH_k^{u,v}$ and $\gamma\in\calH_k^{u+k,v-k+1}$ for $k\geq2$. If there is a
sequence $\alpha_i\in\Omega_{u+1}^{u+v}$ such that
\begin{gather*}
\pi_{u+a,v-a+1}d(\omega+\alpha_i)\lar0\;,\quad 0<a<k\;,\\
\bar Q_{k-2}\pi_{u+k,v-k+1}d(\omega+\alpha_i)\lar0\;,\quad
\Pi_k\pi_{u+k,v-k+1}d(\omega+\alpha_i)\lar\gamma
\end{gather*}
strongly in $\bfOmega$, then $d_k\omega=\gamma$. Moreover, in this case the
sequence $\alpha_i$ can be chosen so that 
\begin{gather*}
\pi_{u+a,v-a+1}d(\omega+\alpha_i)\lar0\;,\quad 0<a<k\;,\\
\pi_{u+k,v-k+1}d(\omega+\alpha_i)\to\gamma
\end{gather*} 
with respect to the \cinf~topology in $\Omega$.
\eprop

The following slightly weaker result will be used as an intermediate step in the proof of
\refp{dk in calHk}.

\lem{dk=0 in calHk}
Let $\gamma\in\calH_k^{u+k,v-k+1}$ for $k\geq2$. If there is some sequence
$\alpha_i\in\Omega_{u+1}^{u+v}$, such that
\begin{gather*}
\pi_{u+a,v-a+1}d\alpha_i\lar0\;,\quad 0<a<k\;,\\
\bar Q_{k-2}\pi_{u+k,v-k+1}d\alpha_i\lar0\;,\quad
\Pi_k\pi_{u+k,v-k+1}d\alpha_i\lar\gamma
\end{gather*}
strongly in $\bfOmega$, then $\gamma=0$.
\elem

Both \refp{dk in calHk} and \refl{dk=0 in calHk} will be proved simultaneously by induction on
$k\geq2$. For the case $k=2$ we need the following.

\lem{Pi2 piu+2,v-2 d tilded1 beta}
We have $\Pi_2\pi_{u+2,v-1}d\tilde d_1\beta=0$ for any $\beta\in\widetilde\calH_1^{u-1,v}$.
\elem

\prf
Write $\beta=\beta'+\beta''$ with
$\beta'\in P(\Omega^{u-1,v})$ and $\beta''\in Q(\Omega^{u,v-1})$.
Then
\begin{align*}
\Pi_2\pi_{u+2,v-2}d\tilde d_1\beta&=
\Pi_2\left(d_{2,-1}\left(d_{1,0}\beta'+d_{0,1}\beta''\right)
+d_{1,0}Q\left(d_{2,-1}\beta'+d_{1,0}\beta''\right)\right)\\
&=\Pi_2\left(\left(d_{2,-1}d_{1,0}+d_{1,0}d_{2,-1}\right)\beta'
+\left(d_{2,-1}d_{0,1}+d_{1,0}^2\right)\beta''\right)\\
&\quad -\Pi_2d_{1,0}\Pi\left(d_{2,-1}\beta'+d_{1,0}\beta''\right)
-\Pi_2d_{1,0}P\left(d_{2,-1}\beta'+d_{1,0}\beta''\right)\\
&=\text{}-\Pi_2d_1\Pi\left(d_{2,-1}\beta'+d_{1,0}\beta''\right)
-\Pi_2Pd_{1,0}P\left(d_{2,-1}\beta'+d_{1,0}\beta''\right)\\
&=0
\end{align*}
by \refe{d0,12=...}, \refl{d10 P=P d10 P, ...} and because $\Pi_2d_{0,1}=\Pi_2d_1=\Pi_2P=0$.
\eprf

\lem{Pi2 piu+2,v-1 d}
Let $\alpha_i$ be a sequence in $\widetilde\calH_1^{u,v}$ such that
$\tilde d_1\alpha_i\to0$ strongly in $\bfOmega$. Then 
$$\Pi_2\pi_{u+2,v-1}d\alpha_i\lar0$$
strongly in $\bfOmega$.
\elem

\prf
Since the image of $\tilde d_1$ is closed and equal to its kernel, the hypothesis implies the
existence of a sequence
$\beta_i\in\widetilde\calH_1^{u-1,v}$ such that $\alpha_i+\tilde d_1\beta_i\to0$ strongly in
$\bfOmega$. On the other hand we have
\begin{align*}
\Pi_2\pi_{u+2,v-1}d&=\Pi_2d_{2,-1}\pi_{u,v}+\Pi_2d_{1,0}\pi_{u+1,v-1}\\
&=\Pi_2d_{2,-1}\pi_{u,v}+\Pi_2\Pi d_{1,0}Q\pi_{u+1,v-1}
\end{align*}
on $\widetilde\calH_1^{u,v}$, and thus the operator
$\Pi_2\pi_{u+2,v-1}d:\widetilde\calH_1^{u,v}\to\calH_2^{u+2,v-1}$ is bounded because
$d_{2,-1}$ and $\Pi d_{1,0}Q$ are bounded operators in $\bfOmega$ by \refl{H1}. Therefore
$$
\Pi_2\pi_{u+2,v-1}d(\alpha_i+\tilde d_1\beta_i)\lar0
$$
strongly in $\bfOmega$. Then the result follows directly from \refl{Pi2 piu+2,v-2 d tilded1 beta}.
\eprf

\prf[Proof of \refl{dk=0 in calHk} for the case $k=2$]
In this case we have $\gamma\in\calH_2^{u+2,v-1}$ and $\alpha_i\in\Omega^{u+1,v-1}$, which satisfy
$$
d_{0,1}\alpha_i\lar0\;,\quad
Qd_{1,0}\alpha_i\lar0\;,\quad
\Pi_2d_{1,0}\alpha_i\lar\gamma
$$
strongly in $\bfOmega$. Since 
$$\Pi d_{1,0}\Pi\alpha_i=d_1\Pi\alpha_i\perp\calH_2^{u+2,v-1}\;,$$
we get $\Pi_2d_{1,0}\widetilde\Pi\alpha_i\to\gamma$ strongly in
$\bfOmega$. But $\Pi_2d_{1,0}P\alpha_i=\Pi_2Pd_{1,0}P\alpha_i=0$ by \refl{d10 P=P d10 P, ...}, and
thus we get 
\begin{equation}\label{e:Pi2 d10 Q alphai lar gamma}
\Pi_2\pi_{u+2,v-1} dQ\alpha_i=\Pi_2 d_{1,0}Q\alpha_i\lar\gamma
\end{equation}
strongly in $\bfOmega$. Now observe that
$Q\alpha_i\in\widetilde\calH_1^{u,v}$, and 
$$
\tilde d_1Q\alpha_i=\widetilde\Pi_{u+1,v}dQ\alpha_i
=d_{0,1}Q\alpha_i+Qd_{1,0}Q\alpha_i
=d_{0,1}\alpha_i+Qd_{1,0}\alpha_i\lar0
$$
because $d_{0,1}Q=d_{0,1}$ and by \refl{d10 P=P d10 P, ...}. Then the
result follows by \refe{Pi2 d10 Q alphai lar gamma} and \refl{Pi2 piu+2,v-1 d}.
\eprf

Now let $\ell\geq2$ and assume that \refl{dk=0 in calHk} holds for $2\leq k\leq\ell$. If $\ell>2$,
assume also that \refp{dk in calHk} holds for $2\leq k<\ell$.

\prf[Proof of \refp{dk in calHk} for $k=\ell$] 
First, we check that the assignment $\omega\mapsto\gamma$, under the conditions in
the statement, defines a map $\calH_\ell^{u,v}\to\calH_\ell^{u+\ell,v-\ell+1}$---observe that, if such a map
is well defined, it is obviously linear---. Suppose there is another
$\gamma'\in\calH_\ell^{u+\ell,v-\ell+1}$ and another sequence
$\alpha'_i\in\Omega_{u+1}^{u+v}$ such that
\begin{gather*}
\pi_{u+a,v-a+1}d\alpha'_i\lar0\;,\quad 0<a<\ell\;,\\
\bar Q_{\ell-2}\pi_{u+\ell,v-\ell+1}d\alpha'_i\lar0\;,\quad
\Pi_\ell\pi_{u+\ell,v-\ell+1}d\alpha'_i\lar\gamma'
\end{gather*}
strongly in $\bfOmega$. Then the sequence $\alpha_i-\alpha'_i\in\Omega_{u+1}^{u+v}$ satisfies
\begin{gather*}
\pi_{u+a,v-a+1}d(\alpha_i-\alpha'_i)\lar0\;,\quad 0<a<\ell\;,\\
\bar Q_{\ell-2}\pi_{u+\ell,v-\ell+1}d(\alpha_i-\alpha'_i)\lar0\;,\quad
\Pi_\ell\pi_{u+\ell,v-\ell+1}d(\alpha_i-\alpha'_i)\lar\gamma-\gamma'
\end{gather*}
strongly in $\bfOmega$. Therefore $\gamma=\gamma'$ by \refl{dk=0 in calHk} for the case $k=\ell$.

Second, we prove that the above map $\calH_\ell^{u,v}\to\calH_\ell^{u+\ell,v-\ell+1}$ is
$d_\ell$; i.e., for each $\omega\in\calH_\ell^{u,v}$, we prove the existence of a sequence
$\alpha_i\in\Omega_{u+1}^{u+v}$ such that
\begin{gather}
\pi_{u+a,v-a+1}d(\omega+\alpha_i)\lar0\;,\quad 0\leq a<\ell\;,
\label{e:piu+a,v-a+1 d(omega+alphai) lar0}\\
\bar Q_{\ell-2}\pi_{u+\ell,v-\ell+1}d(\omega+\alpha_i)\lar0\;,
\label{e:Ql-2 piu+l,v-l+1 d(omega+alphai) lar0}\\
\Pi_\ell\pi_{u+\ell,v-\ell+1}d(\omega+\alpha_i)\lar d_\ell\omega
\in\calH_\ell^{u+\ell,v-\ell+1}\label{e:Pil piu+l,v-l+1 d(omega+alphai) lar0}
\end{gather}
strongly in $\bfOmega$. According to \refe{summarize}, for each
$\omega\in\calH_\ell^{u,v}$ there is a sequence
$\omega_i\in z_\ell^{u,v}$ converging to $\omega$ with respect to the \cinf~topology and such that
$\omega$ and all the $\omega_i$ define the same class $\hat\zeta\in\hat e_{1,\ell}^{u,v}$; thus
all the $\omega_i$ define the same class $\zeta\in e_\ell^{u,v}$. By definition of
$z_\ell^{u,v}$, there is another sequence $\alpha_i\in\Omega^{u+v}_{u+1}$ such that
$\omega_i+\alpha_i\in Z_\ell^{u,v}$. So all the
$\omega_i+\alpha_i$ define the same class
$\xi\in E_\ell^{u,v}$, and the class $d_\ell\xi\in E_\ell^{u+\ell,v-\ell+1}$ is defined by
any of the forms $d(\omega_i+\alpha_i)\in Z_\ell^{u+\ell,v-\ell+1}$. Thus 
$$
\pi_{u+a,v-a+1}d(\omega_i+\alpha_i)=0\;,\quad 0\leq a<\ell\;,
$$
and any of the forms
$$
\pi_{u+\ell,v-\ell+1}d(\omega_i+\alpha_i)\in z_\ell^{u+\ell,v-\ell+1}
$$
define the class $d_\ell\zeta\in e_\ell^{u+\ell,v-\ell+1}$ as well as the class
$d_\ell\hat\zeta\in\hat e_{1,\ell}^{u+\ell,v-\ell+1}$, yielding 
\begin{gather*}
\pi_{u+a,v-a+1}d(\omega_i+\alpha_i)=0\;,\quad 0\leq a<\ell\;,\\
\bar Q_{\ell-2}\pi_{u+\ell,v-\ell+1}d(\omega_i+\alpha_i)=0\;,\\
\Pi_\ell\pi_{u+\ell,v-\ell+1}d(\omega_i+\alpha_i)=d_\ell\omega\in\calH_\ell^{u+\ell,v-\ell-1}
\end{gather*}
independently of $i$. Then~\refe{piu+a,v-a+1
d(omega+alphai) lar0}--\refe{Pil piu+l,v-l+1 d(omega+alphai) lar0} follow by the \cinf~convergence
$\omega_i\to\omega$, as desired. 

Finally we prove the last part of the statement.
Observe that, in fact, the above arguments yield \cinf~convergence in \refe{piu+a,v-a+1
d(omega+alphai) lar0}--\refe{Pil piu+l,v-l+1 d(omega+alphai) lar0}, and also the \cinf~convergence
$$
\bar Q_{\ell-1}\pi_{u+\ell,v-\ell+1}d(\omega+\alpha_i)\lar0\;.
$$
For each $i$, take $\sigma_i^1\in Q_1(\Omega^{u+\ell-1,v-\ell+1})$ satisfying
\begin{equation}\label{e:Pi1 d1,0 sigmai1 = P1 ...}
\Pi_1d_{1,0}\sigma_i^1=P_1\pi_{u+\ell,v-\ell+1}d(\omega+\alpha_i)\;,
\end{equation}
and take $\sigma_i^0\in Q_0(\Omega^{u+\ell,v-\ell})$ such that
\begin{equation}\label{e:d0,1 sigmai0 + P0 d1,0 sigmai1 - P0... lar0}
d_{0,1}\sigma_i^0+P_0d_{1,0}\sigma_i^1
-P_0\pi_{u+\ell,v-\ell+1}d(\omega+\alpha_i)\lar0
\end{equation}
with respect to the \cinf~topology. If $\ell>2$, for each $i$ and $m=2,\ldots,\ell-1$ take some
$\sigma_i^m\in Q_m(\Omega^{u+\ell-m,v-\ell+m})$ such that
$$
d_m\sigma_i^m=P_m\pi_{u+\ell,v-\ell+1}d(\omega+\alpha_i)\;.
$$
By \refp{dk in calHk} for $k<\ell$ there are
sequences $\tau_{i,j}^m\in\Omega_{u+\ell-m+1}^{u+v}$ such that
\begin{gather*}
\pi_{u+\ell-m+a,v-\ell+m-a+1}d(\sigma_i^m+\tau_{i,j}^m)\lar0\;,\quad 0<a<m\;,\\
\pi_{u+\ell,v-\ell+1}d(\sigma_i^m+\tau_{i,j}^m)\lar P_m\pi_{u+\ell,v-\ell+1}d(\omega+\alpha_i)
\end{gather*}
with respect to the \cinf~topology in $\Omega$. Then, for each $i,m$ we can clearly choose $j$
depending on $i,m$ so that $\tau_i^m=\tau_{i,j}^m$ satisfies
\begin{gather}
\pi_{u+\ell-m+a,v-\ell+m-a+1}d(\sigma_i^m+\tau_i^m)\lar0\;,\quad 0<a<m\;,
\label{e:piu+l-m+a,v-l+m-a+1 d(sigmaim + tauim) lar0}\\
\pi_{u+\ell,v-\ell+1}d(\sigma_i^m+\tau_i^m)
-P_m\pi_{u+\ell,v-\ell+1}d(\omega+\alpha_i)\lar0
\label{e:piu+l,v-l+1 d(sigmaim + tauim) - Pm ... lar0}
\end{gather}
with respect to the \cinf~topology. Let 
$$
\beta_i=\alpha_i-\sigma_i^0-\sigma_i^1-\sum_{m=2}^{\ell-1}(\sigma_i^m+\tau_i^m)
\in\Omega_{u+1}^{u+v}\;,
$$
where the last term does not show up if $\ell=2$. 
From~\refe{Pi1 d1,0 sigmai1 = P1 ...}--\refe{piu+l,v-l+1 d(sigmaim + tauim) - Pm ... lar0} we get
\begin{gather*}
\pi_{u+a,v-a+1}d(\omega_i+\beta_i)\lar0\;,\quad 0\leq a<\ell\;,\\
\pi_{u+\ell,v-\ell+1}d(\omega_i+\beta_i)\lar d_\ell\omega
\end{gather*}
with respect to the \cinf~topology in $\Omega$, and the proof is finished. 
\eprf

We already know that both \refp{dk in calHk} and \refl{dk=0 in calHk} hold for $k\leq\ell$, and we
have to prove \refl{dk=0 in calHk} for $k=\ell+1$. The arguments will be similar to the case
$k=2$, and thus we need an appropriate version of \refl{Pi2 piu+2,v-1 d}. In particular,
the generalization of $\widetilde\calH_1^{u,v}$ that fits
our needs turns out to be the following: 
$$
\widetilde\calH_\ell^{u,v}
=P_0(\Omega^{u,v})\oplus\bigoplus_{0<a<\ell}\Omega^{u+a,v-a}
\oplus\bar Q_{\ell-1}(\Omega^{u+\ell,v-\ell})\;.
$$
Let also $\widetilde\calH_\ell^{\cdot,v}=\bigoplus_u\widetilde\calH_\ell^{u,v}$.
We have orthogonal projections $\widetilde\Pi_{\ell;u,v}:\Omega\to\widetilde\calH_\ell^{u,v}$ and
$\widetilde\Pi_{\ell;\cdot,v}:\Omega\lar\widetilde\calH_\ell^{\cdot,v}$ given by
$$
\widetilde\Pi_{\ell;u,v}
=P_0\pi_{u,v}+\sum_{0<a<\ell}\pi_{u+a,v-a}+\bar Q_{\ell-1}\pi_{u+\ell,v-\ell}\;,\quad
\widetilde\Pi_{\ell;\cdot,v}=\sum_u\widetilde\Pi_{\ell;u,v}\;,
$$
and let $\tilde d_\ell=\widetilde\Pi_{\ell;\cdot,v}d:
\widetilde\calH_\ell^{\cdot,v}\to\widetilde\calH_\ell^{\cdot,v}$.

\lem{tilde dl 2=0}
We have $\tilde d_\ell^2=0$
\elem

\prf Consider the following subspaces of $\Omega^{u+v}$:
\begin{gather*}
\calA_\ell^{u,v}=Z_\ell^{u,v}+\overline{B_0^{u,v}}+\Omega_{u+1}^{u+v}
=\left(z_\ell^{u,v}+\overline{b_0^{u,v}}\right)\oplus\Omega_{u+1}^{u+v}\;,\\
\calB^{u,v}=B_0^{u,v}+\Omega_{u+1}^{u+v}=b_0^{u,v}\oplus\Omega_{u+1}^{u+v}\;.
\end{gather*}
First, observe that each $\widetilde\calH_\ell^{u,v}$ is the orthogonal complement of 
$\calA_{\ell-1}^{u+\ell,v-\ell}$ in $\overline{\calB^{u,v}}$, and thus
$\widetilde\calH_\ell^{\cdot,v}$ is the orthogonal complement of
$\calA_{\ell-1}^{\cdot,v-\ell}=\bigoplus_a\calA_\ell^{a,v-\ell}$ in
$\overline{\calB^{\cdot,v}}=\bigoplus_a\overline{\calB^{a,v}}$. So the inclusion
$\widetilde\calH_\ell^{\cdot,v}\hookrightarrow\overline{\calB^{\cdot,v}}$ induces an isomorphism
of topological vector spaces
\begin{equation}\label{e:widetilde calHl.,v=overline calB.,v/calAl-1.,v-l}
\widetilde\calH_\ell^{\cdot,v}\stackrel{\cong}{\lar}
\overline{\calB^{\cdot,v}}/\calA_{\ell-1}^{\cdot,v-\ell}
\end{equation}
whose inverse is induced by the orthogonal projection
$\widetilde\Pi_{\cdot,v}:\overline{\calB^{\cdot,v}}\to\widetilde\calH_\ell^{\cdot,v}$. Second, 
observe that both $\overline{\calB^{u,v}}$ and $\calA_{\ell-1}^{\cdot,v-\ell}$ are subcomplexes of
$(\Omega,d)$. Moreover $\tilde d_\ell$ in $\widetilde\calH_\ell^{\cdot,v}$ clearly corresponds
to the differential map in the quotient complex
$\overline{\calB^{\cdot,v}}/\calA_{\ell-1}^{\cdot,v-\ell}$ via \refe{widetilde calHl.,v=overline
calB.,v/calAl-1.,v-l}, and the result follows.
\eprf

Since $H(\widetilde\calH_1^{\cdot,v},\tilde d_1)=0$, the following two lemmas generalize \refl{Pi2
piu+2,v-2 d tilded1 beta}.

\lem{Pil+1 piu+l,v-l d tildePil d beta}
For any 
$$
\beta\in\bigoplus_{a<\ell}\Omega^{u-1+a,v-a}+\bar Q_{\ell-1}(\Omega^{u-1+\ell,v-\ell})
$$
we have 
$$
\Pi_{\ell+1}\pi_{u+\ell+1,v-\ell}d\widetilde\Pi_{\ell;\cdot,v}d\beta=0\;.
$$
\elem

\prf
By the expression
\begin{multline*}
\bigoplus_{a<\ell}\Omega^{u-1+a,v-a}+\bar Q_{\ell-1}(\Omega^{u-1+\ell,v-\ell})\\
=\bigoplus_{a<\ell}\Omega^{u-1+a,v-a}
+\widetilde\calH_1^{u+\ell-2,v-\ell+1}+(Q_1+\cdots+Q_{\ell-1})(\Omega^{u+\ell-1,v-\ell})
\end{multline*}
it is enough to consider the following three cases. First, assume
$\beta\in\bigoplus_{a<\ell}\Omega^{u-1+a,v-a}$ and let $\beta'=\pi_{u+\ell-2,v-\ell+1}\beta$. We
clearly have
$$
(d-\widetilde\Pi_{\ell;\cdot,v}d)\beta=(\id-\bar Q_{\ell-1})d_{2,-1}\beta'
=(Q_\ell+\Pi_{\ell+1}+\bar P_\ell)d_{2,-1}\beta'\;,
$$
yielding
\begin{align*}
\Pi_{\ell+1}\pi_{u+\ell+1,v-\ell}d\widetilde\Pi_{\ell;\cdot,v}d\beta
&=-\Pi_{\ell+1}\pi_{u+\ell+1,v-\ell}d(Q_\ell+\Pi_{\ell+1}
+\bar P_\ell)d_{2,-1}\beta'\\
&=-\Pi_{\ell+1}d_{1,0}(Q_\ell+\Pi_{\ell+1}+\bar P_\ell)d_{2,-1}\beta'\\
&=-\Pi_{\ell+1}d_1(Q_\ell+\Pi_{\ell+1}+P_1+\cdots+P_\ell)d_{2,-1}\beta'\\
&\quad -\Pi_{\ell+1}d_{1,0}P_0d_{2,-1}\beta'\\
&=0
\end{align*}
by \refl{d10 P=P d10 P, ...}, and because $\Pi_{\ell+1}d_1=0$ and $\Pi_{\ell+1}P_0=0$.

Second, suppose $\beta\in\widetilde\calH_1^{u+\ell-2,v-\ell+1}$ and write $\beta=\beta'+\beta''$
with 
$$
\beta'\in P_0(\Omega^{u+\ell-2,v-\ell+1})\;,\quad\beta''\in Q_0(\Omega^{u+\ell-1,v-\ell})\;.
$$
We clearly have
$$
(\widetilde\Pi_{\ell;\cdot,v}d-\tilde d_1)\beta=(\tilde d_\ell-\tilde
d_1)\beta=(Q_1+\cdots+Q_{\ell-1})(d_{2,-1}\beta'+d_{1,0}\beta'')\;,
$$
yielding
\begin{align*}
\Pi_{\ell+1}\pi_{u+\ell+1,v-\ell}d\widetilde\Pi_{\ell;\cdot,v}d\beta
&=\Pi_{\ell+1}\pi_{u+\ell+1,v-\ell}d\tilde d_1\beta\\
&\quad +\Pi_{\ell+1}d_{1,0}(Q_1+\cdots+Q_{\ell-1})(d_{2,-1}\beta'+d_{1,0}\beta'')\\
&=\Pi_{\ell+1}\Pi_2\pi_{u+\ell+1,v-\ell}d\tilde d_1\beta\\
&\quad +\Pi_{\ell+1}d_1(Q_1+\cdots+Q_{\ell-1})(d_{2,-1}\beta'+d_{1,0}\beta'')\\
&=0
\end{align*}
by \refl{Pi2 piu+2,v-2 d tilded1 beta}.

Third, assume $\beta\in(Q_1+\cdots+Q_{\ell-1})(\Omega^{u+\ell-1,v-\ell})$, which is contained in
$\calH_1^{u+\ell-1,v-\ell}$. Then the result follows
because $\widetilde\Pi_{\ell;\cdot,v}d=\widetilde\Pi_{\ell;\cdot,v}d_1$ on
$\calH_1^{\cdot,v-\ell}$, and
$$
d_1\calH_1^{\cdot,v-\ell}\subset
P_1(\Omega^{\cdot,v-\ell})\perp\widetilde\calH_\ell^{\cdot,v}\;.\qed
$$ 
\renewcommand{\qed}{}\eprf

\lem{Pil+1 piu+l,v-l d alpha}
For $\alpha\in\widetilde\calH_\ell^{u,v}$, if $\tilde d_\ell\alpha=0$, then
$\Pi_{\ell+1}\pi_{u+\ell+1,v-\ell}d\alpha=0$.
\elem

\prf
Write $\alpha=\alpha'+\alpha''+\alpha'''$
with $\alpha'\in P_0(\Omega^{u,v})$, $\alpha''\in\Omega^{u+1,v-1}$ and
$$
\alpha'''\in\bigoplus_{2\leq a<\ell}\Omega^{u+a,v-a}\oplus\bar
Q_{\ell-1}(\Omega^{u+\ell,v-\ell})\;.
$$

Observe that $\alpha'+Q_0\alpha''\in\widetilde\calH_1^{\cdot,v}$. Since $\tilde d_\ell\alpha=0$
and $\tilde d_1=\widetilde\Pi_{1;\cdot,v}\tilde d_\ell$ on $\widetilde\calH_1^{\cdot,v}$, we have
$\tilde d_1(\alpha'+Q_0\alpha'')=0$. Thus there is some $\beta\in\widetilde\calH_1^{u-1,v}$ with
$\tilde d_1\beta=\alpha'+Q_0\alpha''$ because $H(\widetilde\calH_1^{\cdot,v},\tilde d_1)=0$. Then
$\alpha-\tilde d_\ell\beta\in\widetilde\calH_\ell^{u,v}$ satisfies
$$
\pi_{u,v}(\alpha-\tilde d_\ell\beta)=Q_0\pi_{u+1,v-1}(\alpha-\tilde d_\ell\beta)=0\;,
$$
and moreover
$$
\Pi_{\ell+1}\pi_{u+\ell+1,v-\ell}d\alpha
=\Pi_{\ell+1}\pi_{u+\ell+1,v-\ell}d(\alpha-\tilde d_\ell\beta)
$$
by \refl{Pil+1 piu+l,v-l d tildePil d beta}. Therefore we can assume $\alpha'+Q_0\alpha''=0$, and
thus $\alpha'=Q_0\alpha''=0$. With this assumption, it follows that
$\alpha''=(\Pi_1+P_0)\alpha''$ and
$$
d_1\Pi_1\alpha''=\Pi_1d_{1,0}\Pi_1\alpha''=\Pi_1d_{1,0}\alpha''=\Pi_1\pi_{u+2,v-1}d\alpha
=\Pi_1\pi_{u+2,v-1}\tilde d_\ell\alpha=0
$$
by \refl{d10 P=P d10 P, ...}, yielding $Q_1\alpha''=0$.

Take a sequence
$$
\phi_i\in Q_0(\Omega^{u+1,v-2})\subset\widetilde\calH_\ell^{u-1,v}
$$
such that $d_{0,1}\phi_i$ is \cinf~convergent to $P_0\alpha''$.
Then the sequence $\alpha-\tilde d_\ell\phi_i\in\widetilde\calH_\ell^{u,v}$ satisfies
\begin{align*}
\Pi_{\ell+1}\pi_{u+\ell+1,v-\ell}d\alpha
&=\Pi_{\ell+1}\pi_{u+\ell+1,v-\ell}d(\alpha-\tilde d_\ell\phi_i)\\
&\lar\Pi_{\ell+1}\pi_{u+\ell+1,v-\ell}d(\Pi_1\alpha''+\alpha''')
\end{align*}
by \refl{Pil+1 piu+l,v-l d tildePil d beta}. So
$$
\Pi_{\ell+1}\pi_{u+\ell+1,v-\ell}d\alpha
=\Pi_{\ell+1}\pi_{u+\ell+1,v-\ell}d(\Pi_1\alpha''+\alpha''')\;,
$$
and thus we can also assume $P_0\alpha''=0$.

For each $k=1,\dots,\ell-1$ there is some
$\sigma^k\in Q_k(\Omega^{u-k+1,v+k-2})$
with $d_k\sigma^k=P_k\alpha''$. As above, from the existence of such a $\sigma^1$ we
can assume $P_1\alpha''=0$ by \refl{Pil+1 piu+l,v-l d tildePil d beta} since
$d_1=\pi_{u,v-1}\Pi_1\tilde d_\ell$ on
$\calH_1^{u,v-1}$. If $\ell>2$, by \refp{dk in calHk} for $k=2,\ldots,\ell-1$ there is a sequence
$\tau^k_i\in\Omega^{u+v-1}_{u-m+2}$ such that
\begin{gather*}
\pi_{u-k+a+1,v+k-a-1}d(\sigma^k+\tau^k_i)\lar0\;,\quad 0<a<k\;,\\
\pi_{u+1,v-1}d(\sigma^k+\tau^k_i)\lar P_k\alpha''
\end{gather*}
with respect to the \cinf~topology in $\Omega$. We can thus suppose $P_k\alpha''=0$ for
such a $k$ because 
$$
\Pi_{\ell+1}\pi_{u+\ell+1,v-\ell}d\widetilde\Pi_{\ell;\cdot,v}d(\sigma^k+\tau^k_i)=0
$$
by \refl{Pil+1 piu+l,v-l d tildePil d beta}. Therefore 
\begin{equation}\label{e:alpha'' in ...}
\alpha''\in\calH_\ell^{u+1,v-1}\oplus\bigoplus_{k=2}^{\ell-1}Q_k(\Omega^{u+1,v-1})\;,
\end{equation}
where the last term does not show up if $\ell=2$.

Now the condition $\tilde d_\ell\alpha=0$ can be written as
\begin{equation}\label{e:d(alpha''+alpha''')}
\pi_{u+1+a,v-a}d(\alpha''+\alpha''')
=\bar Q_{\ell-1}\pi_{u+\ell+1,v-\ell}d(\alpha''+\alpha''')=0\;,\quad0<a<\ell\;.
\end{equation}
Observe that~\refe{d(alpha''+alpha''')} summarizes the conditions of the first part of \refp{dk in calHk} for
$k=2,\dots,\ell$, with $\omega=\alpha''$, the constant sequence $\alpha_i=\alpha'''$, and 
$\gamma=0$ if $2\leq k<\ell$. Since $\alpha''\in\calH_2^{u+1,v-1}$ by~\refe{alpha'' in ...}, we
get inductively on $k=2,\dots,\ell-1$ that $\alpha''\in\calH_k^{u+1,v-1}$ and $d_k\alpha''=0$
by~\refe{alpha'' in ...},~\refe{d(alpha''+alpha''')} and \refp{dk in calHk}. Hence
$\alpha''\in\calH_\ell^{u+1,v-1}$ by~\refe{alpha'' in ...},
and thus
$$
\Pi_\ell\pi_{u+\ell+1,v-\ell}d\alpha=\Pi_\ell\pi_{u+\ell+1,v-\ell}d(\alpha''+\alpha''')
=d_\ell\alpha''\perp\calH_{\ell+1}
$$
by~\refe{d(alpha''+alpha''')} and \refp{dk in calHk} for $k=\ell$, and the result follows.
\eprf

We also need the following Hodge theory for the complex $(\widetilde\calH_\ell^{\cdot,v},\tilde
d_\ell)$. Let $\tilde\delta_\ell=\widetilde\Pi_{\ell;\cdot,v}\delta$ on
$\widetilde\calH_\ell^{\cdot,v}$, and set $\widetilde D_\ell=\tilde
d_\ell+\tilde\delta_\ell$ and
$\widetilde\D_\ell=\widetilde D_\ell^2=\tilde\delta_\ell\tilde
d_\ell+\tilde d_\ell\tilde\delta_\ell$. Such a
$\tilde\delta_\ell$ is adjoint of
$\tilde d_\ell$ in $\widetilde\calH_\ell^{\cdot,v}$ with respect to the $L^2$~inner product, and
thus $\widetilde D_\ell$ and $\widetilde\D_\ell$ are symmetric unbounded operators in the
$L^2$~completion $L^2\widetilde\calH_\ell^{\cdot,v}$. 

\lem{tildeDl is self-adjoint}
The operator $\widetilde D_\ell$ is essentially self-adjoint in
$L^2\widetilde\calH_\ell^{\cdot,v}$.
\elem

\prf 
By Theorem~2.2 in \cite{Chernoff}, $D=d+\delta$ is essentially
self-adjoint in $\bfOmega$. Then, by using e.g. Lemma~XII.1.6--(c) in
\cite{DunfordSchwartz}, so is $\widetilde\Pi_{\ell;\cdot,v} D\widetilde\Pi_{\ell;\cdot,v}$
because $\widetilde\Pi_{\ell;\cdot,v}$ is a bounded self-adjoint
operator on $\bfOmega$. But $\widetilde\Pi_{\ell;\cdot,v} D\widetilde\Pi_{\ell;\cdot,v}$ is equal
to $\widetilde D_\ell$ in $L^2\widetilde\calH_\ell^{\cdot,v}$ and vanishes in its orthogonal
complement. Hence $\widetilde D_\ell$ is essentially self-adjoint.
\eprf

\lem{WrtildecalHl}
$D\widetilde\Pi_{\ell;\cdot,v}-\widetilde\Pi_{\ell;\cdot,v}D\widetilde\Pi_{\ell;\cdot,v}$ defines
a bounded operator on $\bfOmega$.
\elem

\prf
We have
\begin{align*}
D\widetilde\Pi_{\ell;\cdot,v}-\widetilde\Pi_{\ell;\cdot,v}D\widetilde\Pi_{\ell;\cdot,v}
&=\widetilde P_0(\delta_{-1,0}+\delta_{-2,1}\pi_{\cdot,v-1})
+\delta_{-2,1}\pi_{\cdot,v}\\
&\quad +(\Pi_\ell+\bar
P_{\ell-1})(d_{1,0}\pi_{\cdot,v-\ell}+d_{2,-1}\pi_{\cdot,v-\ell+1})\\ 
&\quad +d_{2,-1}\pi_{\cdot,v-\ell}\;.
\end{align*}
But 
$$
\widetilde P_0\delta_{-1,0}\pi_{\cdot,v}=\widetilde P_0\delta_{-1,0} P_0\pi_{\cdot,v}\;,\quad
P_0d_{1,0}\pi_{\cdot,v-\ell}=P_0d_{1,0}\widetilde P_0\pi_{\cdot,v-\ell}
$$
on $\widetilde\calH_\ell^{\cdot,v}$. Then the result follows by \refl{H1}-(i).
\eprf

For each positive integer $r$, define the norm $\|\cdot\|'_r$ 
on $\widetilde\calH_\ell^{\cdot,v}$ by setting
$$\|\phi\|'_r=\left\|(\id+\widetilde D_\ell)^r\phi\right\|\;,$$
and let $W^k\widetilde\calH_\ell^{\cdot,v}$ be the
corresponding completion of $\widetilde\calH_\ell^{\cdot,v}$. Then the
following result follows directly from~\refl{WrtildecalHl}.

\cor{WrtildecalHl}
The restriction of each $r$th Sobolev norm 
$\|\cdot\|_r$ to $\widetilde\calH_\ell^{\cdot,v}$ is equivalent 
to the norm $\|\cdot\|'_r$. Thus $W^k\widetilde\calH_\ell^{\cdot,v}$ is the
closure of $\widetilde\calH_\ell^{\cdot,v}$ in $W^k\Omega$.
\ecor

\cor{Hodge theory for tildecalHl}
The Hilbert space $L^2\widetilde\calH_\ell^{\cdot,v}$ has a complete
orthonormal system $\{\phi_i : i=1,2,\ldots\}\subset\widetilde\calH_\ell^{\cdot,v}$,
consisting of eigenvectors of $\widetilde\D_\ell$, so that
the corresponding eigenvalues satisfy
$0\leq\lambda_1\leq\lambda_2\leq\cdots$ with $\lambda_i\uparrow\infty$
if $\dim\widetilde\calH_\ell^{\cdot,v}=\infty$; thus all of these eigenvalues have finite
multiplicity. We also have the orthogonal decomposition
$$
\widetilde\calH_\ell^{\cdot,v}=(\ker\tilde d_\ell\cap\ker\tilde\delta_\ell)\oplus
\im\tilde d_\ell\oplus\im\tilde\delta_\ell\;,
$$
with
\begin{align*}
\ker\widetilde\D_\ell&=\ker\tilde d_\ell\cap\ker\tilde\delta_\ell\;,\\
\im\widetilde\D_\ell&=\im\tilde d_\ell\oplus\im\tilde\delta_\ell\;,\\
\ker d_\ell&=(\ker\tilde d_\ell\cap\ker\tilde\delta_\ell)\oplus\im\tilde d_\ell\;,\\
\ker\delta_\ell&=(\ker\tilde d_\ell\cap\ker\tilde\delta_\ell)\oplus\im\tilde\delta_\ell\;.
\end{align*}
\ecor

\prf
\refc{WrtildecalHl} implies that each inclusion $W^{r+1}\widetilde\calH_\ell^{\cdot,v}\hookrightarrow
W^r\widetilde\calH_\ell^{\cdot,v}$ is a compact operator, and
$\bigcap_rW^r\widetilde\calH_\ell^{\cdot,v}=\widetilde\calH_\ell^{\cdot,v}$.  Then the result follows by
Proposition~2.44 in \cite{AlvTond1} and
\refl{tildeDl is self-adjoint}.
\eprf

Contrary to the case of $(\widetilde\calH_1^{\cdot,v},\tilde d_1)$, it may easily happen
that the complex
$(\widetilde\calH_\ell^{\cdot,v},\tilde d_\ell)$ has non-trivial cohomology. But we still can 
finish the proof of \refl{dk=0 in calHk}.

\prf[Proof of \refl{dk=0 in calHk} for the case $k=\ell+1$]
Observe that the strong convergence
\begin{gather*}
\pi_{u+a,v-a+1}d\alpha_i\lar0\;,\quad 0<a\leq\ell\;,\\
\bar Q_{\ell-1}\pi_{u+\ell+1,v-\ell}d\alpha_i\lar0
\end{gather*}
in $\bfOmega$ just means the strong convergence
$\tilde d_\ell\widetilde\Pi_{\ell;u,v}\alpha_i\to0$. Write
$\widetilde\Pi_{\ell;u,v}\alpha_i=\phi_i+\psi_i$ with $\phi_i\in\ker\tilde d_\ell$ and
$\psi_i\in\im\tilde\delta_\ell$, according to \refc{Hodge theory for tildecalHl}. Then $\tilde
d_\ell\psi_i\to0$ strongly in $\bfOmega$ by \refl{tilde dl 2=0}, yielding $\psi_i\to0$ strongly in
$\bfOmega$ by \refc{Hodge theory for tildecalHl}. Moreover
$$
\Pi_{\ell+1}\pi_{u+\ell+1,v-\ell}d\alpha_i
=\Pi_{\ell+1}\pi_{u+\ell+1,v-\ell}d\widetilde\Pi_{\ell;u,v}\alpha_i
=\Pi_{\ell+1}\pi_{u+\ell+1,v-\ell}d\psi_i\lar0
$$
by \refl{Pil+1 piu+l,v-l d alpha} and because any linear map
$\widetilde\calH_\ell^{\cdot,v}\to\calH_{\ell+1}^{\cdot,v-\ell}$ is continuous with respect to
the $L^2$~norms since $\calH_{\ell+1}$ is of finite dimension. Therefore
$\gamma=0$ as desired.
\eprf

This finishes the proof of \refp{dk in calHk}, which has the following consequence.

\cor{deltak in calHk}
Let $\omega\in\calH_k^{u,v}$ and $\gamma\in\calH_k^{u-k,v+k-1}$ for $k\geq2$. If there is a
sequence $\beta_i\in\bigoplus_{a>0}\Omega^{u-a,v+a}$ such that
\begin{gather*}
\pi_{u-a,v+a-1}\delta(\omega+\beta_i)\lar0\;,\quad 0<a<k\;,\\
\bar P_{k-2}\pi_{u-k,v+k-1}\delta(\omega+\beta_i)\lar0\;,\quad
\Pi_k\pi_{u-k,v+k-1}\delta(\omega+\beta_i)\lar\gamma
\end{gather*}
strongly in $\bfOmega$, then $\delta_k\omega=\gamma$. Moreover, in this case the
sequence $\beta_i$ can be chosen so that 
\begin{gather*}
\pi_{u-a,v+a-1}\delta(\omega+\beta_i)\lar0\;,\quad 0<a<k\;,\\
\pi_{u-k,v+k-1}\delta(\omega+\beta_i)\lar\gamma
\end{gather*} 
with respect to the \cinf~topology in $\Omega$.	
\ecor

\prf
We can assume that $M$ is oriented by using the two fold covering of orientations with standard
arguments. Then it is easy to check that the Hodge star operator, $\star:\Omega\to\Omega$,
satisfies $\star\calH_k=\calH_k$, and $\star d_k=(-1)^{r+1}\delta_k\star$ on $\calH_k^r$ for each
integer $r$. Then the result follows from \refp{dk in calHk}.
\eprf

Now \reft{calHk} follows directly from \refp{dk in calHk} and \refc{deltak in calHk} by induction
on $k$.

\ssec{rescaling}{Estimates of the rescaled Laplacian}

The rescaled Laplacian $\D_h$ is the square of the
``rescaled Dirac operator''
$D_h=d_h+\delta_h$, which will be used here too.
The sum of \refe{dh} and \refe{deltah} gives
\begin{equation}\label{e:D h}
D_h=D_0+hD_\perp+h^2F\;,
\end{equation}
where 
$$
D_0=d_{0,1}+\delta_{0,-1}\;,\quad
D_\perp=d_{1,0}+\delta_{-1,0}\;,\quad
F=d_{2,-1}+\delta_{-2,1}\;,
$$
Let also $\D_\perp=D_\perp^2$.

\begin{Lem}[{\'Alvarez-Kordyukov \cite[Remark 3.5]{AlvKordy1}}]\label{l:Dperp D0+D0 Dperp}
There is a zero order differential operator $B$
on $\Omega$ such that
$$D_\perp D_0+D_0D_\perp=BD_0+D_0B^\ast\;.$$
\end{Lem}

\prop{Delta h}
There is some $C>0$ such that\footnote{Recall that, for self-adjoint operators $A,B$ in a Hilbert space
$H$, the inequality $A\leq B$ is defined in the sense of quadratic forms: $\langle Au,u\rangle\leq\langle
Bu,u\rangle$ for all $u\in H$.}
$$
\D_h\geq \frac{1}{2}\D_0+\frac{1}{2}h^2\D_\perp-Ch^2
$$
for $h$ small enough.
\eprop

\prf Consider the operators $B,F$ given by \refl{Dperp D0+D0 Dperp} and \refe{D h}.
Since $B,F$ are of order zero, there is
some $C'>0$ such that
$B^\ast B,F^2\leq C'$. Because $D_0$ is symmetric, we get
\begin{align*}
h\,|\langle(BD_0+D_0B^\ast)\omega,\omega\rangle|
&\leq 2h\,|\langle D_0\omega,B\omega\rangle|\\
&\leq 2h\,\|D_0\omega\|\,\|B\omega\|\\
&\leq\frac{1}{4}\,\|D_0\omega\|^2+4h^2\,\|B\omega\|^2\\
&=2\left\langle\left(\frac{1}{4}\D_0+4h^2B^\ast B\right)\omega,\omega\right\rangle
\end{align*}
for all $\omega\in\Omega$, yielding
$$h\,|BD_0+D_0B^\ast|\leq\D_0+h^2B^\ast B\leq\D_0+C'h^2.$$
Similarly we get
$$\begin{array}{c}
\left|FD_0+D_0F\right|\leq\D_0+F^2\leq\D_0+C'\;,\\[6pt]
\left|FD_\perp+D_\perp F\right|\leq\D_\perp+F^2\leq\D_0+C'\;.
\end{array}$$
Therefore, from \refe{D h} and \refl{Dperp D0+D0 Dperp} we get
\begin{align*}
\D_h&=\D_0+h^2\D_\perp+h^4F^2+h(BD_0+D_0B^\ast)\\
&\phantom{=}\text{}+h^2(D_0F+FD_0)+h^3(D_\perp F+FD_\perp)\\
&\geq\D_0+h^2\D_\perp+h^4C'-\frac{1}{4}\D_0-C'h^2\\
&\phantom{=}\text{}-h^2\left(\D_0+C'\right)-h^3\left(\D_\perp+C'\right)\\
&\geq\frac{1}{2}\D_0+\frac{1}{2}h^2\D_\perp-Ch^2
\end{align*}
for some $C>0$ and all $h$ small enough.
\eprf

\begin{proof}[Proof of \refmt{asymptotics of eigenforms}] In the case $k=1$, 
\refe{asymptotics of eigenforms}
just means
$\left\langle\D_{h_i}\omega_i,\omega_i\right\rangle\to0$.
Therefore
$$
\left\langle\left(\frac{1}{2}\D_0+\frac{1}{2}h_i^2\D_\perp-Ch_i^2\right)\omega_i,\omega_i\right\rangle\lar0
$$
by \refp{Delta h}. Hence
\begin{equation}\label{e:<Delta0 omegai,omegai> to 0}
\langle\D_0\omega_i,\omega_i\rangle\lar0
\end{equation}
and $\langle\D_\perp\omega_i,\omega_i\rangle$ is uniformly bounded since both
$\D_0$ and $\D_\perp$ are positive operators. It follows that
$\omega_i$ is uniformly bounded in $W^1\Omega$. Therefore some subsequence of
$\omega_i$ is weakly convergent in $W^1\Omega$ (and thus strongly convergent in
$\bfOmega$) to some $\omega\in W^1\Omega$. From \refe{<Delta0 omegai,omegai> to 0} we also
get that $\|D_0\omega_i\|\to0$. So 
$D_0\omega_i\to0$ strongly in $\bfOmega$, yielding
$\omega\in\ker\bfD_0$ because $\bfD_0$ is a closed operator in $\bfOmega$. But $\ker\bfD_0=L^2\calH_1$ by
\refe{L2calH1}. Thus the result follows for $k=1$.

For $k=2$, it follows from \refe{asymptotics of eigenforms} that 
$$\left\|d_{h_i}\omega_i\right\|\in o(h_i)\;,\quad
\left\|\delta_{h_i}\omega_i\right\|\in o(h_i)\;,$$
yielding that 
$$
\left(\frac{1}{h_i}d_{0,1}+d_{1,0}+h_id_{2,-1}\right)\omega_i
\lar0\;,\quad
\left(\frac{1}{h_i}\delta_{0,-1}+\delta_{-1,0}+h_i\delta_{-2,1}\right)\omega_i
\lar0\;,
$$
strongly in $\bfOmega$ by \refe{dh} and \refe{deltah}. Hence
$$
\Pi\left(d_{1,0}+h_id_{2,-1}\right)\omega_i
\lar0\;,\quad
\Pi\left(\delta_{-1,0}+h_i\delta_{-2,1}\right)\omega_i
\lar0
$$ 
strongly in $\bfOmega$ as well, and thus so does the sequence
$\Pi D_\perp\omega_i$. Then 
$$
D_1\Pi\omega_i=\Pi D_\perp\Pi\omega_i
=\Pi D_\perp\omega_i-\Pi D_\perp\widetilde\Pi\omega_i\lar0
$$
strongly in $\bfOmega$ by \refl{H1}-(i). It follows that $\omega\in\ker\bfD_1$
because $\bfD_1$ is a closed operator in $L^2\calH_1$. But $\ker\bfD_1=\calH_2$ by \reft{H1}-(iii), and the
result follows for $k=2$.

For the case $k>2$, we can assume $\omega_i\in\Omega^r$ and 
$\omega\in\calH_2^{u,v}$ for some integers $u+v=r$. Let
$\omega_i^a=\pi_{a,r-a}\omega_i$ for each integer $a$, and set
$$
\omega_i'=\sum_{a\geq0}h_i^{-a}\omega_i^{u+a}\;,\quad
\omega_i''=\sum_{a\geq0}h_i^{-a}\omega_i^{u-a}\;.
$$
Now, by \reft{calHk}, the result follows from the following claim.

\begin{claim}\label{claim:omega'}
For $0<a<k$, we have 
$$\pi_{u+a,v-a+1}d\omega_i'\lar0\;,\quad
\pi_{u-a,v+a-1}\delta\omega_i''\lar0\;,$$
strongly in $\bfOmega$.
\end{claim}

Clearly 
$$\pi_{u,v+1}d\omega_i'=d_{0,1}\omega_i^u\;,\quad
\pi_{u,v-1}\delta\omega_i''=\delta_{0,-1}\omega_i^u\;.$$
Thus both of these components converge strongly to zero because $\omega\in L^2\calH_1$. 

To prove Claim~\ref{claim:omega'} for other bihomogeneous components observe
that, again from \refe{asymptotics of eigenforms}, both 
$\left\|d_{h_i}\omega_i\right\|$ and $\left\|\delta_{h_i}\omega_i\right\|$ are
in $o\left(h_i^{k-1}\right)$.  Then
\begin{equation}\label{e:omega'}
\left\|h_i^2d_{2,-1}\omega_i^{b-2}+h_id_{1,0}\omega_i^{b-1}
+d_{0,1}\omega_i^b\right\|
\in o\left(h_i^{k-1}\right)\;,
\end{equation}
\begin{equation}\label{e:omega''}
\left\|h_i^2\delta_{-2,1}\omega_i^{b+2}+h_i\delta_{-1,0}\omega_i^{b+1}
+\delta_{0,-1}\omega_i^b\right\|
\in o\left(h_i^{k-1}\right)\;,
\end{equation}
for every integer $b$, by considering bihomogeneous components of 
$d_{h_i}\omega_i$ and $\delta_{h_i}\omega_i$. Now 
\begin{align*}
\pi_{u+1,v}d\omega_i'&=d_{1,0}\omega_i^u+h_i^{-1}d_{0,1}\omega_i^{u+1}\;,\\
\pi_{u-1,v}\delta\omega_i''&=
\delta_{-1,0}\omega_i^u+h_i^{-1}\delta_{0,-1}\omega_i^{u-1}\;.
\end{align*}
Both of these components strongly converge to zero in $\bfOmega$ too by 
\refe{omega'} and \refe{omega''}, since so does
$h_id_{2,-1}\omega_i^{u-1}$ and $h_i\delta_{-2,1}\omega_i^{u+1}$ because
$d_{2,-1}$ and $\delta_{-2,1}$ are of order zero and $\|\omega_i\|=1$. 

The other bihomogeneous components of $d\omega_i'$ and $\delta\omega_i''$ are the following ones,
where $a\geq2$,
\begin{align*}
\pi_{u+a,v-a+1}d\omega_i'&=h_i^{-a+2}d_{2,-1}\omega_i^{u+a-2}+h_i^{-a+1}d_{1,0}\omega_i^{u+a-1}
+h_i^{-a}d_{0,1}\omega_i^{u+a}\;,\\
\pi_{u-a,v+a-1}\delta\omega_i''&=
h_i^{-a+2}\delta_{-2,1}\omega_i^{u-a+2}+h_i^{-a+1}\delta_{-1,0}\omega_i^{u-a+1}
+h_i^{-a}\delta_{0,-1}\omega_i^{u-a}\;,
\end{align*}
which strongly converge to zero in $\bfOmega$ for $a<k$ by
\refe{omega'} and \refe{omega''}. This finishes the proof of
Claim~\ref{claim:omega'}.\eprf

\begin{proof}[Proof of \refmt{small eigenvalues}] 
First, we can assume the metric is bundle-like by \refe{N'rh}. So we can apply the results
of this section.

If we had a strict inequality ``$<$'' in~\refe{small eigenvalues k} for some
$k\geq2$, by the isomorphism $\calH_k^r\cong E_k^r$ there are sequences $\omega_i\in\Omega^r$ and
$h_i\downarrow0$ such that $\|\omega_i\|=1$, $\omega_i\perp\calH_k$, and
$$\langle\D_{h_i}\omega_i,\omega_i\rangle\in O(h_i^{2k})\;.$$
But then we get a contradiction by \refmt{asymptotics of eigenforms}. So inequality ``$\geq$''
holds in~\refe{small eigenvalues k} for all $k\geq2$.

The proof of ``$\geq$'' in \refe{small eigenvalues 1}  follows with the same arguments since
$\widehat E_1^r\cong\calH_1^r$, which is of finite dimension if and only if so is $ L^2\calH_1^r$.

For $k\geq2$, inequality ``$\leq$'' of
\refe{small eigenvalues k} in \refmt{small eigenvalues} follows directly from 
\refc{bfEkr Hausdorff of finite dim} and \reft{Ek=bfEk}, as was pointed out in \refr{leq}.

Now observe that, for each $h>0$ and each $\omega\in\calH_1^r$, we have 
$$D_h\omega=hD_\perp\omega+h^2F\omega\;,$$
according to \refe{D h}. Therefore the inequality ``$\leq$'' in \refe{small eigenvalues 1} 
follows from the isomorphism $\calH_1^r\cong\widehat E_1^r$ by using the well known variational formula
$N_h^r(\lambda)=\sup_V\dim V$, where $V$ runs over the subspaces of $\Omega^r$ satisfying
$$\langle\D_h\omega,\omega\rangle\leq\lambda\,\|\omega\|^2$$
for all $\omega\in V$.
\eprf

\sec{Forman}{Forman's nested sequence}

This section is devoted to the proof of \refmt{frakHk0,p}. Thus let $\calF$ be a Riemannian foliation of
dimension $p$ on a closed manifold $M$. We need the following characterization of ${\mathfrak H}_2$, which
is weaker than \refe{frakHk} for $k=2$.

\begin{claim}\label{claim:frakH2} 
A differential form $\omega\in\Omega$ is in ${\mathfrak H}_2$ if and only if it has extensions
$\tilde\omega_1(h),\tilde\omega_2(h)\in\Omega[h]$ satisfying
\begin{equation}\label{e:frakH2}
d_h\tilde\omega_1(h)\in h^2\Omega[h]\;,\quad\delta_h\tilde\omega_2(h)\in h^2\Omega[h]\;.
\end{equation}
\end{claim}

According to \refe{frakHk}, it is enough to prove the ``if'' part of Claim~\ref{claim:frakH2}. We can assume 
$$
\tilde\omega_1(h)=\omega+h\omega_1\;,\quad\tilde\omega_2(h)=\omega+h\omega_2
$$
for some $\omega_1,\omega_2\in\Omega$ because $d_h(h^2\Omega[h])$ and $\delta_h(h^2\Omega[h])$ are contained
in $h^2\Omega[h]$. On the other hand, since ${\mathfrak H}_2$ is a bigraded subspace of $\Omega$, we can
suppose $\omega\in\Omega^{u,v}$ for some $u,v$. Then it easily follows from
\refe{frakH2} that
$\omega_1\in\Omega^{u+1,v-1}$ and $\omega_2\in\Omega^{u-1,v+1}$. Furthermore
we can assume $\delta_{0,-1}\omega_1=d_{0,1}\omega_2=0$ by \reft{leafwise Hodge}. Hence the extension
$$
\tilde\omega(h)=\omega+h(\omega_1+\omega_2)
$$
of $\omega$ is easily seen to satisfy \refe{frakHk} for $k=2$, and thus $\omega\in{\mathfrak H}_2$, finishing
the proof of Claim~\ref{claim:frakH2}

The statement of Claim~\ref{claim:frakH2} seems to hold also for ${\mathfrak H}_k$ with
$k>2$, but the proof can not be so easy. 

By \refmt{small eigenvalues} and \refe{frakH1=calH1, ...}, we have ${\mathfrak H}_2^{0,p}=\calH_2^{0,p}=0$
if $E_2^{0,p}=0$. Therefore we can assume $E_2^{0,p}\neq0$ to prove \refmt{frakHk0,p}. 
According to \cite{Masa92} and \cite{Alv7}, this assumption implies that \calF\ is orientable and
$E_2^{0,p}\cong\R$. So $\calH_2^{0,p}\cong\R$ by \refmt{small eigenvalues}, and thus either ${\mathfrak
H}_2^{0,p}(g)=0$ or ${\mathfrak H}_2^{0,p}(g)=\calH_2^{0,p}(g)$ by \refe{frakH1=calH1, ...}. 

Recall from \cite{Rummler} that the {\em characteristic form\/}, determined by \calF\ and a metric $g$ on
$M$, is the unique differential form $\chi\in\Omega^{0,p}$ whose restriction to the leaves is the leafwise
volume form. If $g$ is a bundle-like metric, then $\delta_{0,-1}$ corresponds
to the leafwise coderivative by restriction to the leaves \cite{Alv7}, \cite{AlvKordy1}, yielding
$\delta_{0,-1}\chi=0$, and thus $\chi\in{\mathfrak H}_1^{0,p}(g)$.

To prove \refmt{frakHk0,p}-(i) just choose the bundle-like metric $g$ so that $d_{1,0}\chi=0$, which can be
done by using Sullivan's purification \cite{Sullivan79} (see also
\cite{Masa92} and \cite{Alv7}). Hence $\chi\in{\mathfrak H}_2^{0,p}(g)$ by Claim~\ref{claim:frakH2}, yielding
${\mathfrak H}_2^{0,p}(g)\neq0$. 

To prove \refmt{frakHk0,p}-(ii), let us begin with a bundle-like metric $g$ satisfying 
\refmt{frakHk0,p}-(i), and the corresponding bigrading of $\Omega$ and decomposition of $d$ and $\delta$ as
sum of bihomogeneous components. The hypothesis
$\bar0_1^{0,p}\neq0$ means that $d_{0,1}\Omega^{0,p-1}$ is not closed in $\Omega^{0,p}$, and thus we can
take some
$\alpha\in\overline{d_{0,1}\Omega^{0,p-1}}\setminus d_{0,1}\Omega^{0,p-1}$. Take also some $\epsilon>0$
small enough so that
$\chi+\epsilon\alpha=f\chi$ for some positive function $f$. Therefore  $\chi'=f\chi$ is the characteristic
form of some bundle-like metric $g'$ on $M$. Such a $g'$ can be chosen to define the same bigrading on
$\Omega$ as
$g$,  yielding the same decomposition of $d$ as sum of bihomogeneous components. We have $\chi'\in{\mathfrak
H}_1^{0,p}(g')=\calH_1^{0,p}(g')$. Moreover, since $\alpha$ defines a non-trivial class 
$$
[\alpha]\in\overline{d_{0,1}\Omega^{0,p-1}}/d_{0,1}\Omega^{0,p-1}=\bar o_1^{0,p}\cong\bar 0_1^{0,p}
$$
and since $H^0(\bar o_1^{\cdot,p})=H^0(\bar 0_1^{\cdot,p})=0$ by \reft{H1}-(vi), we get 
$$0\neq d_1[\alpha]=[d_{1,0}\alpha]\in\bar o_1^{1,p}\cong\bar 0_1^{1,p}\;.$$
So
$$
d_{1,0}\chi'=d_{1,0}(\chi+\epsilon\alpha)=\epsilon d_{1,0}\alpha\in\overline{d_{0,1}\Omega^{1,0}}\setminus
d_{0,1}\Omega^{1,0}\;,
$$
yielding $\chi'\in\calH_2^{0,p}\setminus{\mathfrak H}_2^{0,p}(g')$. Therefore
${\mathfrak H}_2^{0,p}(g')\neq\calH_2^{0,p}(g')$, and thus ${\mathfrak H}_2^{0,p}(g')=0$.

\begin{ack}
The second author gratefully acknowledges the hospitality and support of the University of Santiago de
Compostela.
\end{ack}

\end{document}